\documentclass[letterpaper]{article}
\usepackage[english]{babel}
\usepackage{amssymb,amsmath}
\usepackage{bbm}     
\usepackage{ifpdf}
\usepackage{ifvtex}
\ifvtexpdf\pdftrue\fi
\ifpdf
\usepackage[pdftex]{graphicx}
\usepackage[pdftex,pdftitle={On the Geometric Dilation of Closed Curves, Graphs, and Point Sets},bookmarks]{hyperref}
\else
\usepackage[dvips]{graphicx}
\fi

\makeatletter
\def\ps@headings{\let\@mkboth\markboth
  \def\@oddfoot{\hfil \thepage \hfil}
  \def\@evenfoot{\hfil \thepage \hfil}
  \def\@evenhead{\em Preprint submitted to Elsevier Science\hfil \today}
  \def\@oddhead{\em Preprint submitted to Elsevier Science\hfil \today}
}
\makeatother
\newtheorem{thm}{Theorem}
\newtheorem{lem}{Lemma}
\newtheorem{cor}{Corollary}
\newtheorem{claim}{Claim}
\newenvironment{pf}{\begin{par}\noindent{\bf Proof.}}{\end{par}\par\medskip}

\providecommand{\qed}{\Box}

\newcommand{\abs}[1]{\left| #1 \right|}
\newcommand{\edge}[1]{#1} 
\newcommand{\pseudofrac}[2]{#1 / #2}

\newcommand{\eps}[0]{\varepsilon}
\newcommand{\realset}[0]{\mathbbm{R}}

\newcommand{\rationalset}[0]{\mathbbm{Q}}
\newcommand{\C}{\mathcal{C}}
\newcommand{\M}{\mathcal{M}}
\newcommand{\D}{\mathcal{D}}
\newcommand{\halfC}[0]%
{{\textstyle\frac{\abs{C}}{2}}}
\newcommand{\halfCtext}[0]%
{|C|/2}
\newcommand{\intd}[0]%
{\mathrm{d}}
\newcommand{\DeltaBound}[0]%
{\underline{\Delta}}
\newcommand{\vect}{\genfrac(){0pt}{}}
\newcommand{\comment}[1]{}

\begin{document}


\title{On the Geometric Dilation of Closed Curves, Graphs, and Point Sets%
\footnotemark[1]}

\author{Adrian Dumitrescu\footnotemark[2]
\and
Annette Ebbers-Baumann\footnotemark[3]
\and
Ansgar Gr{\"u}ne\footnotemark[3]
\footnotemark[5]
\and
Rolf Klein\footnotemark[3]
\footnotemark[6]
\and
G{\"u}nter Rote\footnotemark[4]
}
\date{August 25, 2005}
\maketitle
\renewcommand{\thefootnote}{\fnsymbol{footnote}}
\footnotetext[1]{Some results of this article were presented at
  the 21st European Workshop on Computational Geometry (EWCG '05)\cite{dgr-ilbgd-05},
 others at the 9th Workshop on Algorithms and Data Structures (WADS '05)\cite{degkr-gdhc-05}.}
\footnotetext[2]{Computer Science, University of Wisconsin--Milwaukee,
    3200 N. Cramer Street, Milwaukee, WI 53211, USA; {\tt ad@cs.uwm.edu}}
\footnotetext[3]{Institut f{\"u}r Informatik I, Universit{\"a}t Bonn, R{\"o}merstra{\ss}e 164, D - 53117 Bonn, Germany;\\
  {\{\texttt{ebbers}, \texttt{gruene}, \texttt{rolf.klein}\}\tt @cs.uni-bonn.de}}
\footnotetext[4]{Freie Universit{\"a}t Berlin, Institut f{\"u}r Informatik,
  Ta\-ku\-stra\ss e 9, D-14195 Berlin, Germany;\\
  {\tt rote@inf.fu-berlin.de}}
\footnotetext[5]{Ansgar Gr{\"u}ne was partially supported by a DAAD PhD-grant.}
\footnotetext[6]{Rolf Klein was partially supported by DFG-grant KL 655/14-1.}
\renewcommand{\thefootnote}{\arabic{footnote}}

\begin{abstract}
  Let $G$ be an embedded planar graph whose edges are curves.
  The \emph{detour} between two points $p$ and $q$ (on edges or vertices)
  of~$G$ is the ratio between the length of a shortest path connecting~$p$
  and~$q$ in~$G$ and their Euclidean distance~$|pq|$.
  The maximum detour over all pairs of points is called
  the \emph{geometric dilation} $\delta(G)$.

  Ebbers-Baumann, Gr{\"u}ne and Klein have shown that every finite point set is
  contained in a planar graph whose geometric dilation is at most $1.678$,
  and some point sets require graphs with dilation $\delta\ge \pi/2
  \approx 1.57$.
  They conjectured that the lower bound is not tight.

  We use new ideas like the halving pair transformation,
  a disk packing result and
  arguments from convex geometry,
  to prove this conjecture.
  The lower bound is improved to $(1+ 10^{-11})\pi/2$.
  The proof relies on halving pairs,
  pairs of points dividing a given closed curve~$C$ in two parts
  of equal length, and their minimum and maximum distances~$h$ and~$H$.
  Additionally, we analyze curves of constant halving distance ($h=H$),
  examine the relation of~$h$ to other geometric quantities
  and prove some new dilation bounds.

  \vspace{1em}
  \noindent{\em Key words:}
  computational geometry, convex geometry, convex curves, dilation, distortion,
  detour, lower bound, halving chord, halving pair, Zindler curves
\end{abstract}


\section{Introduction}
\begin{figure}[htb]
  \centering
\noindent
\vbox{\halign{\hfil#\hfil\cr
  \includegraphics{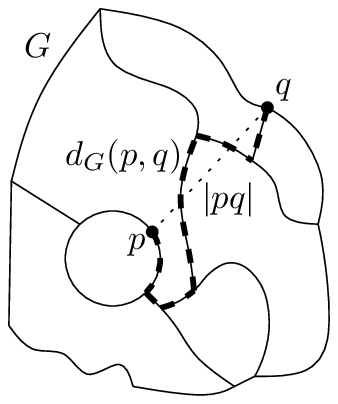}\cr
(a)\cr}}
\hfil
\vbox{\halign{\hfil#\hfil\cr
  \includegraphics{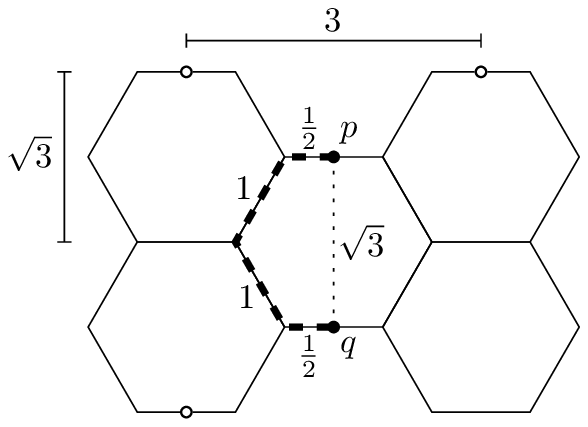}
\cr
(b)\cr}}

  \caption{
    (a) The shortest path (dashed) between $p$ and $q$ in the graph $G$ and
      the direct distance $|pq|$.
    (b) Three points (drawn as empty circles) embedded in a hexagonal grid
        with geometric dilation $\sqrt{3}\approx 1.732$.}
  \label{DilationExample-fig}
\end{figure}

Consider a planar graph $G$ embedded in $\realset^2$, whose edges are
curves\footnote{For simplicity we assume here that the curves are
piecewise continuously differentiable, but most of the
proofs can be extended to arbitrary rectifiable curves.}
that do not intersect.
Such graphs arise naturally in the study
of transportation networks, like waterways, railroads or streets.
For two points, $p$ and $q$ (on edges or vertices) of $G$, the
\emph{detour} between $p$ and $q$ in $G$ is defined as
\[
  \delta_G(p,q)
  =
  \frac{d_G(p,q)}{|pq|}
\]
where $d_G(p,q)$ is the shortest path length in~$G$ between~$p$ and~$q$,
and $|pq|$ denotes the Euclidean distance,
see Figure~\ref{DilationExample-fig}a for an illustration.
Good transportation networks should
have small detour values.
In a railroad system, access is only possible at stations,
the vertices of the graph.
Hence, to measure the quality of such networks, we can take the maximum
detour over all pairs of vertices.
This results in the well-known concept of {\em graph-theoretic dilation}
studied extensively in the literature on spanners,
see \cite{e-sts-00} for a survey.

However, if we consider a system of urban streets,
houses are usually spread everywhere along the streets.
Hence, we have to take into account not only the vertices of the graph
but all the points on its edges.
The resulting supremum value is the {\em geometric dilation}
\begin{equation*}
  \delta(G)
  :=
  \sup\limits_{p,q \in G} \delta_G(p,q)
  =
  \sup\limits_{p,q \in G} \frac{d_G(p,q)}{|pq|}
\end{equation*}
on which we concentrate in this article.
Several papers \cite{ekll-faadp-04,lms-cmdsr-02,akks-cdpc-02}
have shown how to efficiently compute the geometric dilation of
polygonal curves.
Besides this the geometric dilation was studied in
differential geometry and knot theory
under the notion of distortion,
see e.g.~\cite{glp-smvr-81,ks-dtk-98}.

\begin{figure}[htbp]
  \begin{center}
    \includegraphics{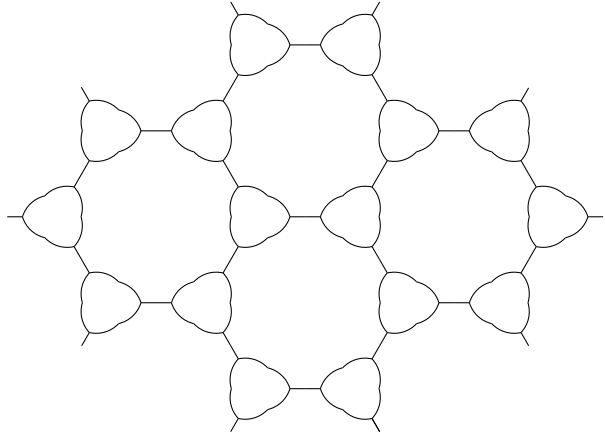}
  \end{center}
  \caption{\label{ImprovedHexGrid-fig}%
    A section of the grid from \cite{egk-gdfps-03} which has small dilation, less than 1.678.}
\end{figure}
Ebbers-Baumann et al.~\cite{egk-gdfps-03} recently
considered the problem of constructing
a graph of lowest possible geometric dilation
containing a given finite point set on its edges.
Even for three given points this is not a trivial task.
For some examples, clearly a Steiner-tree with three straight line segments is
optimal.
In other cases a path consisting of straight and curved pieces is better,
but it is not easy to prove its optimality.

Therefore, Ebbers-Baumann et al.\ concentrated
on examining the dilation necessary to embed
any finite point set, i.e. the value
\begin{equation*}
  \Delta
  :=
  \sup\limits_{P \subset \realset^2,\;P\;\mathrm{finite}}
  \inf\limits_{G \supset P,\;G\;\mathrm{finite}} \delta(G)\;.
\end{equation*}
The infimum is taken over all embedded planar graphs~$G$
with a finite number of vertices, where the edges may be
curves like discussed above.
For example, a scaled hexagonal grid
can clearly be used to embed any finite subset of
$\rationalset \times \sqrt{3}\rationalset$
like the example in Figure~\ref{DilationExample-fig}b.
The geometric dilation $\sqrt{3}\approx 1.732$
of this grid is attained between two midpoints
of opposite edges of a hexagon.

Ebbers-Baumann et al.\ introduced the improved grid
shown in Figure~\ref{ImprovedHexGrid-fig} and proved that
its dilation is less than~$1.678$.
They showed that a slightly perturbed version
of the grid can be used to embed any finite point set.
Thereby they proved $\Delta<1.678$.

They also derived that
$\Delta \geq \pi/2$, by showing that a graph
$G$ has to contain a cycle to embed a certain
point set $P_5$ with low dilation, and by using
that the dilation of every closed curve\footnote{In this paper
we use the
notions ``cycle'' and ``closed curve'' synonymously.}~$C$
 is bounded by $\delta(C)\geq\pi/2$.

They conjectured that this lower bound is not tight.
It is known that circles are the only
cycles of dilation $\pi/2$,
see \cite[Corollary 23]{egk-gdcpc-04ii}, \cite[Corollary 3.3]{acfgh-cmmke-02},
\cite{ks-dtk-98}, \cite{glp-smvr-81}.
And intuition suggests that
one cannot embed complicated point sets with dilation $\pi/2$
because every face of the graph would have to be a circular disk.
This idea would have to be formalized
and still does not rule out that every point
set could be embedded with dilation arbitrarily close to $\pi/2$.

Therefore, we need the result presented in Section~\ref{RingResult-sec}.
We show that cycles with dilation
\emph{close} to $\pi/2$ are \emph{close} to circles, in some well-defined
sense (Lemma~\ref{Ring-lem}).
The lemma can be seen as an instance of a \emph{stability result}
for the geometric inequality $\delta(C)\geq\pi/2$,
 see~\cite{g-sgi-93} for a survey. Such
results complement geometric inequalities
 (like the isoperimetric inequality between the area
and the perimeter of a planar region)
with statements of the following
kind: When the inequality is fulfilled ``almost'' as an equation, the object
under investigation is ``close'' to the object or class of objects for which
the inequality is tight.
An important idea in the proof of this stability result
is the decomposition of any closed curve $C$
into the two cycles $C^*$ and $M$ defined in Section~\ref{MAndCStar-sec}.

In Section~\ref{MainResult-sec} we use Lemma~\ref{Ring-lem}
to relate the dilation problem to a certain problem of
packing and covering the plane by disks.
By this we prove our main result $\Delta \geq (1+10^{-11})\pi/2$.
The proof also relies on the notion of halving pairs
and their distance, the halving distance,
introduced by Ebbers-Baumann et al.~\cite{egk-gdcpc-04ii}
to facilitate the dilation analysis of closed curves.

In Section~\ref{ConstantHalvingDistance-sec} we analyze curves
of constant halving distance,
an analog to the well-known curves of constant width.
Understanding curves of constant halving distance and their properties
is a key point in designing networks with small geometric dilation.
For example, the grid structure in Figure~\ref{DilationExample-fig}b
is constructed by replicating such a curve at each vertex of a regular
hexagonal grid, thereby improving the dilation
from~$\sqrt{3}\approx 1.732$ to~$1.678$.
Curves of constant halving distance were already discovered in
1921~\cite{z-ukg-21}%
\footnote{We would like to thank Salvador Segura-Gomis for pointing
  this out.}, and, as we will explain, they are related to other
interesting geometric notions, such as curves of constant width and
Stanis\l aw Ulam's Floating Body Problem.

From the viewpoint of convex geometry it is interesting to
consider the relations of the minimum and maximum halving distance,
$h$ and~$H$, to other
geometric quantities of a given convex closed curve $C$.
In Section~\ref{hAndOther-sec} we give first results
in this direction.

In Section~\ref{DilationBounds-sec} we use some of them
to derive a new upper bound on the geometric dilation
of closed convex curves.
We also prove several new dilation bounds for polygons.

\section{Basic Definitions and Properties}
\label{BasicDefinitions-sec}%
An important special case of the planar graphs embedded in $\realset^2$
are simple%
\footnote{A curve is called {\em simple} if it has
no self-intersections.}
{\em closed curves} or {\em cycles} for short.
Let $C$ be such a closed curve.
By~$|C|$ we denote its length.
Shortest path distance $d_C(p,q)$,
detour $\delta_C(p,q)$ and geometric dilation $\delta(C)$
are defined like in the case of arbitrary graphs.
Often we will use
a bijective arc-length-parameterization $c: [0,|C|) \to C$.
This implies $d_C(c(s),c(t)) = \min\left(|t-s|,|C|-|t-s|\right)$,
and $|\dot{c}(t)|=1$ wherever the derivative exists.

\begin{figure}[htbp]
  \begin{center}
    \includegraphics{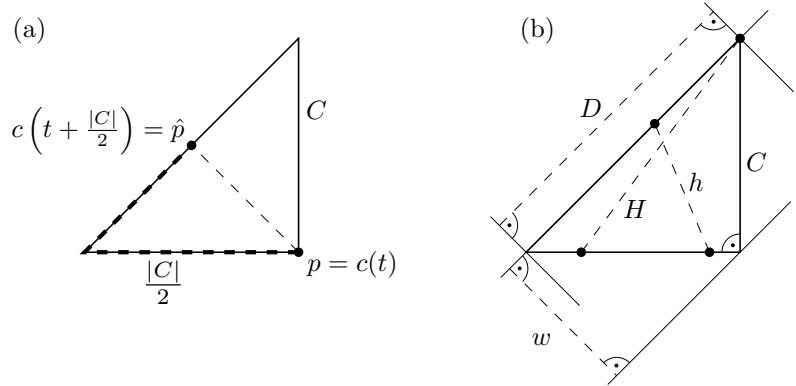}
  \end{center}
  \caption{\label{HalvingDistance-fig}%
    (a)~A halving chord $p\hat{p}$.
    (b)~Diameter $D$, width $w$, minimum and maximum halving distance
    $h$ and $H$ of an isosceles, right-angled triangle.}%
\end{figure}

Consider Figure~\ref{HalvingDistance-fig}.
Two points $p=c(t)$ and $\hat{p}=c(t+\halfCtext)$ on $C$ that divide
the length of $C$ in two equal parts form a \emph{halving pair} of $C$.
Here and later on, $t+\halfCtext$ is calculated modulo $|C|$.
The segment which connects $p$ and $\hat{p}$ is a \emph{halving chord},
and its length is the corresponding \emph{halving distance}.
We write $h=h(C)$ and $H=H(C)$ for the \emph{minimum} and \emph{maximum
halving distance} of~$C$.
Furthermore, we will consider the {\em diameter} $D := \max\{|pq|, p,q\in C\}$
of a closed curve $C$
and the {\em width} $w$ of a convex cycle $C$ which is the minimum
distance of two parallel lines enclosing $C$.

The main link between dilation and halving distance is the inequality
\begin{equation}\label{InequalityDilationh-equ}
  \delta(C)
  \geq
  \frac{|C|}{2h}
\end{equation}
which follows immediately from the definitions
because the right-hand side equals the detour of
a halving pair of minimal distance~$h$.
If $C$ is convex, equality is attained.
\begin{lem}\label{DilationConvex-lem}{\rm \cite[Lemma 11]{egk-gdcpc-04ii}}
If $C$ is a closed convex curve, its dilation $\delta(C)$ is attained by a
halving pair, i.e. $\delta(C) = |C|/2h$.
\end{lem}
%

\section{Midpoint Curve and Halving Pair Transformation}
\label{MAndCStar-sec}%
\begin{figure}[ht]
  \centering
  \includegraphics{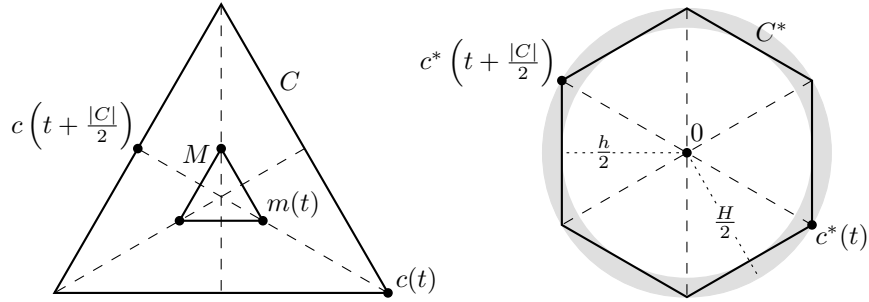}
  \caption{\label{ThreeCurves-fig} An equilateral triangle $C$,
    and the two curves $C^*$ and $M$ derived from $C$.}
\end{figure}

Let $C$ be a closed curve and let $c(t)$ be an arc-length parameterization.
The two curves derived from $C$ illustrated in Figure~\ref{ThreeCurves-fig}
turn out to be very useful for both,
proving the stability result in Section~\ref{RingResult-sec}
and analyzing the curves of constant halving distance
in Section~\ref{ConstantHalvingDistance-sec}.
The \emph{midpoint curve} $M$ is
formed by the midpoints of the halving chords of $C$,
and is given by the parameterization
\begin{equation}
  \label{M-equ}
  m(t) := \frac{1}{2} \left( c(t) + c(t+\halfC) \right)
  .
\end{equation}
The second curve $C^*$ is the result of the {\em halving pair transformation}
introduced in~\cite{egk-gdcpc-04ii}.
It is obtained by translating all the halving chords so that their midpoints
are located in the origin. Then, $C^*$ is the curve described by the endpoints
of these translated chords. This results in the parameterization
\begin{equation}
  \label{CStar-equ}
  c^*(t) := \frac{1}{2} \left( c(t) - c(t+\halfC) \right)
  .
\end{equation}
Note that $c^*(t)$ is half the vector connecting the
corresponding halving pair.
By definition, ${c^*(t)=-c^*(t+\halfCtext)}$, hence $C^*$ is centrally
symmetric.
The curve $C^*$ has the same set of halving distances as $C$; thus,
$h(C^*)=h(C)=h$ and $H(C^*)=H(C)=H$.
Furthermore, it is contained in an $(H/h)$-ring, see Figure~\ref{ThreeCurves-fig}.
An {\em $\eta$-ring} is the closed region between two concentric circles
where the outer radius equals $\eta$ times the inner radius.

The parameterization of the midpoint curve satisfies $m(t)=m(t+\halfCtext)$, and thus,
$M$ traverses the same curve twice when $C$ and $C^*$ are traversed once.
We define $|M|$ as the length of the curve $m(t)$ corresponding to one
traversal, i.e., the
parameter interval is $[0, |C|/2]$.

The halving pair transformation decomposes the curve
$C$ into two components, from which $C$ can be reconstructed:
\begin{equation}
  \label{Decomposition-equ}
  c(t) = m(t) + c^*(t)
  , \qquad
  c\left(t+\halfC\right) = m(t) - c^*(t)
\end{equation}
This is analogous to the decomposition of a function into an even
and an odd function, or writing a matrix as a sum of a symmetric and a
skew-symmetric matrix.

A key fact in the proof of our main result,
Theorem~\ref{ImprovedLowerBound-thm},
is the following lemma, which we think
is of independent interest.
It provides an upper bound on the length
$\abs{M}$ of the midpoint curve in terms of $\abs{C}$ and $\abs{C^*}$.
\begin{lem} \label{LengthInequalityThreeCurves-lem}
$
  \;\;\;
  4|M|^2 + |C^*|^2 \leq |C|^2.
$
\end{lem}

\begin{pf}
As mentioned in the introduction, we assume that $C$ is piecewise
continuously differentiable.
However, the following proof can be extended to arbitrary rectifiable curves
(i.e., curves with finite
length).

Using the linearity of the scalar product
and $\abs{\dot{c}(t)}=1$, we obtain
\begin{eqnarray}
  \label{MCStarOrthogonal-equ}
  \langle \dot m(t), \dot c^*(t) \rangle
  & \stackrel{\textrm{\scriptsize \thetag{\ref{M-equ}},\thetag{\ref{CStar-equ}}}}{=} &
  \frac{1}{4} \left\langle \dot c(t)+\dot c\left(t+\halfC\right),
               \dot c(t)-\dot c\left(t+\halfC\right) \right\rangle
  \\ \nonumber & = &
  \frac{1}{4} \left (|\dot c(t)|^2 -\left|\dot c\left(t+\halfC\right)\right|^2 \right )
  =
  \frac{1}{4} (1-1)=0
  .
\end{eqnarray}
This means that the derivative vectors $\dot c^*(t)$ and $\dot m(t)$
are always orthogonal, thus
\thetag{\ref{Decomposition-equ}} yields
\begin{equation*}
  \abs{\dot{m}(t)}^2 + \abs{\dot c^*(t)}^2
  =
  \abs{\dot{c}(t)}^2
  =
  1.
\end{equation*}
This implies
\begin{eqnarray}
  \nonumber
  |C|& = &
  \int_0^{|C|} \sqrt{\abs{\dot{m}(t)}^2 +
  \abs{\dot c^*(t)}^2} \,\intd t
  \\ \label{LengthCDecomposition-equ} & \geq &
  \sqrt
  {
    \left(\int_0^{|C|} \abs{\dot{m}(t)} \intd t\right)^2
    +
    \left(\int_0^{|C|} \abs{\dot c^*(t)} \intd t\right)^2
  }
  =
  \sqrt{4\abs{M}^2 + |C^{*}|^2}
\end{eqnarray}
The above inequality --- from which the lemma follows --- can be seen by
a geometric argument: the left integral
\[
  \int_0^{|C|} \sqrt{\abs{\dot{m}(t)}^2 +
  \abs{\dot c^*(t)}^2} \,\intd t
\]
  is the length of the curve
\[
  \gamma(s) := \left( \int_0^s \abs{\dot{m}(t)} \intd t,
  \int_0^s \abs{\dot c^*(t)} \intd t\right),
\]
  while the right expression
\[
  \sqrt{
    \left(\int_0^{|C|} \abs{\dot{m}(t)} \intd t\right)^2
    +
    \left(\int_0^{|C|} \abs{\dot c^*(t)} \intd t\right)^2
  }
\]
equals the distance of its end-points
$\gamma(0)=(0,0)$ and $\gamma\left(|C|\right)$.
\hfill$\qed$
\end{pf}

\begin{cor} \label{LengthInequalityCCStar-cor}
$\;\;\; |C^*| \leq |C|$.
\end{cor}

\section{Stability Result for Closed Curves}\label{RingResult-sec}
In this section, we prove that a simple closed curve $C$
of low dilation (close to $\pi/2$) is close to being a circle.
To this end, we first show that $C^*$ is close to a circle,
i.e. $H/h$ is close to 1.
Then, we prove that the length of the midpoint curve is small.
Combining both statements yields the desired result.

We use the following lemma to find an upper bound on the ratio~$H/h$.
It extends an inequality of Ebbers-Baumann et al.~\cite[Theorem 22]{egk-gdcpc-04ii}
to non-convex cycles.
\begin{lem} \label{LowerDilationBoundC-lem}
The geometric dilation $\delta(C)$
of any closed curve $C$ satisfies
\begin{equation*}
\delta(C) \geq \arcsin{\frac{h}{H}} +
\sqrt { \left (\frac {H}{h} \right )^2 -1 } .
\end{equation*}
This bound is tight.
\end{lem}
\begin{figure}[htbp]
  \begin{center}
    \includegraphics{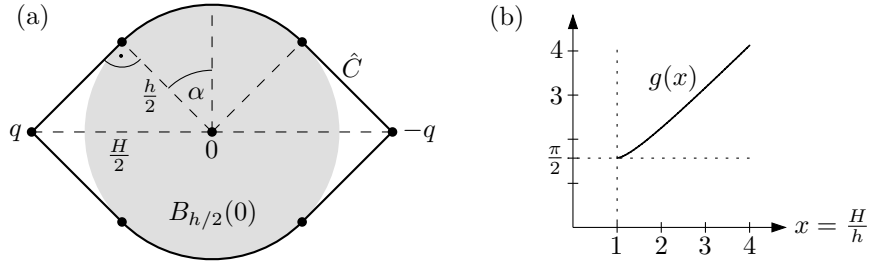}
  \end{center}
  \caption{\label{CapCurveAndBound-fig}%
    (a)~The curve~$\hat{C}$ is the shortest curve enclosing~$B_{\frac{h}{2}}(0)$
        and connecting~$q$ and~$-q$.
    (b)~A plot of the lower bound on geometric dilation depending on~$x=H/h$.}%
\end{figure}
The function~$g(x)=\arcsin{1/x}+\sqrt{x^2-1}$ on the right side is displayed in Figure~\ref{CapCurveAndBound-fig}b.
It starts from $g(1)=\pi/2$ and is increasing on $[1,\infty)$.
This shows that a cycle $C$ of small dilation
has a ratio~$H/h$ close to 1.
\begin{pf}
{From} the definition of dilation and halving pairs
and Corollary~\ref{LengthInequalityCCStar-cor}, we get
\[
  \delta(C)
  \stackrel{\textrm{\scriptsize \thetag{\ref{InequalityDilationh-equ}}}}{\geq}
  \frac{|C|}{2h}
  \stackrel{\textrm{\scriptsize Cor.~\ref{LengthInequalityCCStar-cor}}}{\geq}
  \frac{|C^*|}{2h}
  .
\]
We have seen above that
$h(C^*)=h(C)=h$ and $H(C^*)=H(C)=H$.
Then, by the symmetry of $C^*$, any two points $p$ and $-p$ on $C^*$ form
a halving pair of distance $2\abs{p} \geq h$.
Hence, $C^*$ contains the disk $B_{h/2}(0)$
of radius $h/2$ around the origin.
On the other hand, $C^*$ has to connect some halving pair $(q,-q)$
of maximum distance~$H$.
Therefore, $C^*$ has at least the length of the curve $\hat{C}$
depicted in Figure~\ref{CapCurveAndBound-fig}.
Basic trigonometry yields
\[
  |C^*|
  \geq
  |\hat{C}|
  =
  2h \left ( \arcsin {\frac{h}{H}} +
  \sqrt { \left (\frac {H}{h} \right )^2 -1 } \right ),
\]
which completes the proof.
Kubota~\cite{k-uee-23} used similar arguments to prove
that the length of any convex closed curve~$C$
of width~$w$ and diameter~$D$ satisfies
$\abs{C} \geq 2w\arcsin(w/D) + 2\sqrt{D^2-w^2}$.
Here and in the dilation bound,
equality is attained by $\hat{C}$.
\hfill$\qed$
\end{pf}
In order to prove that a closed curve of small dilation
lies in a thin ring, we still need an upper bound on the
length $|M|$ of the midpoint cycle.
We achieve it by using Lemma~\ref{LengthInequalityThreeCurves-lem}.
\begin{lem}%
  \label{BoundLengthM-lem}
  If $\delta(C) \leq (1+\eps)\pseudofrac{\pi}{2}$,
  then $|M| \leq (\pi h/2)\sqrt{2\eps + \eps^2}$.
\end{lem}
\begin{pf}
  By the assumption and because the dilation of $C$
  is at least the detour of a halving pair
  attaining minimum distance~$h$,
  we have
  $
    (1+\eps)\pi/2
    \geq
    \delta(C)
    \stackrel{\textrm{\scriptsize \thetag{\ref{InequalityDilationh-equ}}}}{\geq}
    |C|/2h
  $, implying
  \begin{equation}
    \label{UpperBoundLengthC-equ}
    |C| \leq (1 + \eps) \pi h .
  \end{equation}
  As seen before, $C^*$ encircles but does not enter
  the open disk~$B_{h/2}(0)$
  of radius~$h/2$ centered at the origin~$0$.
  It follows that the length of~$C^*$ is at least
  the perimeter of~$B_{h/2}(0)$:
  \begin{equation}
    \label{LowerBoundLengthCStar-equ}
    |C^*| \geq \pi h.
  \end{equation}
  By plugging everything together, we get
  \begin{equation*}
    |M|
    \stackrel{\mbox{\scriptsize Lemma~\ref{LengthInequalityThreeCurves-lem}}}{\leq}
    \frac{1}{2}\sqrt{|C|^2-|C^*|^2}
    \stackrel{\mbox{\scriptsize \thetag{\ref{UpperBoundLengthC-equ}},\thetag{\ref{LowerBoundLengthCStar-equ}}}}{\leq}
    \frac{1}{2}\pi h\sqrt{(1+\eps)^2-1}
    =
    \frac{\pi h}{2} \sqrt{2\eps+\eps^2},
  \end{equation*}
  which concludes the proof of Lemma~\ref{BoundLengthM-lem}.
\hfill$\qed$
\end{pf}
Intuitively it is clear (remember the definitions (\ref{M-equ}), (\ref{CStar-equ}) and Figure~\ref{ThreeCurves-fig})
that the upper bound on $H/h$ from Lemma~\ref{LowerDilationBoundC-lem}
and the upper bound on $|M|$ of Lemma~\ref{BoundLengthM-lem}
imply that the curve~$C$ is contained in a thin ring,
if its dilation is close to $\frac{\pi}{2}$.
This is the idea behind
Lemma~\ref{Ring-lem},
the main result of this section.
To prove it, we will apply the following well-known fact, see e.g.~\cite{sa-ics-00}, to $M$.
\begin{lem}\label{CurveInDisk-lem}
 Every closed curve~$C$ can be enclosed in
 a circle of radius $|C|/4$.
\end{lem}
\begin{pf}
Fix a halving pair $(p,\hat{p})$ of $C$.
Then by definition, for any $q \in C$, we have
$ |pq|+|q\hat{p}| \leq d_C(p,q) + d_C(q,\hat{p}) = d_C(p,\hat{p}) = |C|/2 $.
It follows that $C$ is contained in an ellipse with foci $p$ and $\hat{p}$
and major axis $|C|/2$.
This ellipse is included in a circle with
radius $|C|/4$, and the lemma follows.
\hfill$\qed$
\end{pf}
\begin{lem} \label{Ring-lem}
Let $C \subset \realset^2$ be any simple closed curve
with dilation $\delta(C) \leq (1 + \eps)\pi/2$ for $\eps \leq 0.0001$.
Then $C$ can be enclosed in a $(1+3\sqrt{\eps})$-ring.
This bound cannot be improved apart from the coefficient of $\sqrt{\eps}$.
\end{lem}
\begin{pf}
By Lemma~\ref{CurveInDisk-lem}, the midpoint cycle $M$ can
be enclosed in a circle of radius $|M|/4$ and some center~$z$.
By the triangle inequality,
we immediately obtain
\begin{eqnarray*}
    \abs{c(t)-z}
    & \stackrel{\mbox{\scriptsize (\ref{Decomposition-equ})}}{=} &
    \abs{m(t)+c^*(t)-z}
    \leq
    \abs{c^*(t)} + \abs{m(t)-z}
    \leq \frac{H}{2} + \frac{|M|}{4},
    \\
    \abs{c(t)-z}
    & \stackrel{\mbox{\scriptsize (\ref{Decomposition-equ})}}{=} &
    \abs{m(t)+c^*(t)-z}
    \geq
    \abs{c^*(t)}-\abs{m(t)-z}
    \geq \frac{h}{2} - \frac{|M|}{4}.
\end{eqnarray*}
Thus, $C$ can be enclosed
in the ring between two concentric circles with radii
$R=\pseudofrac{H}{2} + \pseudofrac{|M|}{4}$ and
$r=\pseudofrac{h}{2} - \pseudofrac{|M|}{4}$ centered at~$z$.
To finish the proof, we have to bound the ratio~$R/r$.
For simplicity, we prove only the asymptotic bound
$R/r\le 1+O(\sqrt{\eps})$.
The proof of the precise bound, which includes all numerical estimates, is given in
Appendix~\ref{PreciseBound-sec}.

Assume $H/h=(1+\beta)$.
Lemma~\ref{LowerDilationBoundC-lem} and power series expansion
yield the approximate lower bound
\begin{equation}\label{AsymptoticBound-equ}
  \delta(C)
  \geq
  \arcsin \frac{1}{1+\beta} + \sqrt{(1+\beta)^2-1}
  = \frac\pi 2 + \frac{2\sqrt2}3 \beta^{3/2} - O(\beta^{5/2}).
\end{equation}
With our initial assumption
$\delta(C) \leq \frac{\pi}{2}(1+\eps)$
we get therefore
$\beta = O( \eps^{2/3})$.
Lemma~\ref{BoundLengthM-lem} implies
$ |M| = O(h\sqrt{\eps})$
which yields
\[
\frac{R}{r}
=\frac{H/2 +|M|/4}{h/2 -|M|/4}
=\frac{h(1+\beta) + |M|/2}{h - |M|/2}
   \leq \frac{1+ O(\eps^{2/3}) + O(\eps^{1/2})}{1- O(\eps^{1/2})}
=   1+O( \eps^{1/2})
,
\]
completing the proof of the asymptotic bound in Lemma \ref{Ring-lem}.
\hfill$\qed$
\end{pf}
The lemma can be extended to a larger, more practical range of
$\eps$, by increasing the coefficient of~$\sqrt{\eps}$.

{\bf Tightness of the bound in Lemma \ref{Ring-lem}.}
\label{Tightness-sec}
\begin{figure}[htbp]
  \centerline{\includegraphics{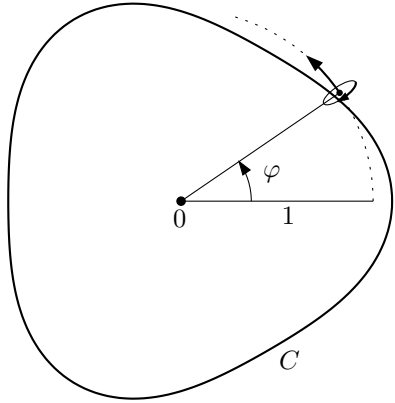}}
  \caption{A moon's orbit; the figure shows the curve $C$ for $s=0.1$.}
  \label{TightExample-fig}
\end{figure}
The curve $C$ defined by the parameterization $c(\varphi)$ below
and illustrated in Figure~\ref{TightExample-fig}
shows that the order of magnitude in the bound cannot be improved.
Note that $c(\varphi)$ is not an arc-length parameterization.
\[
  c(\varphi) := \binom{\cos \varphi}{\sin \varphi} (1+s \cos{3 \varphi})
  + \binom{-\sin \varphi}{\cos \varphi} (-\frac{s}{3} \sin {3 \varphi}).
\]
This curve is the path of a moon moving around the earth on a small
elliptic orbit with major axis~$2s$ (collinear to the
line earth--sun), with a frequency three times that
of the earth's own circular orbit around the sun.

Here, we only sketch how to bound the dilation and
halving distances of this curve.
The details are given in Appendix~\ref{DetailedTightnessExample-sec}.
A ring with outer radius $R$ and inner radius $r$
containing the curve
satisfies $R/r \geq (1+s)/(1-s)=1+\Theta(s)$.
One can show that the length is
bounded by $|C| \leq 2 \pi + O(s^2)$,
and the halving distances are bounded by
$2-O(s^3)\le h<H \leq 2+O(s^3)$.
If $s$ is not too large, $C$ is convex, and
this implies by Lemma~\ref{DilationConvex-lem} that
the dilation is given by
\[
  \delta(C)
  =
  \frac {|C|/2}h \leq \frac{\pi+O(s^2)}{2-O(s^3)}
  =
  (1+O(s^2)) \frac{\pi}{2}
  .
\]
Thus, we have dilation
$\delta=(1+\eps)\pi/2$ with $\eps=O(s^2)$,
but
the ratio of the radii of the enclosing ring is
$1+\Theta(s)=1+\Omega(\sqrt\eps)$.
A more careful estimate shows that this ratio is
$1+\frac{3}{2}\sqrt\eps+O(\eps)$.
Thus, the coefficient~3 in Lemma~\ref{Ring-lem}
cannot be improved very much.

\section{Improved Lower Bound on the Dilation of Finite Point Sets}
\label{MainResult-sec}%
We now combine Lemma~\ref{Ring-lem} with a disk packing result
to achieve the desired new lower bound on $\Delta$.
A (finite or infinite) set $\D$ of disks in the plane with disjoint
interiors is called a \emph{packing}.
\begin{thm} \label{DiskPacking-thm}  {\rm (Kuperberg, Kuperberg,
Matou\v sek and Valtr~\cite{kkmv-atpe-99})}
Let $\D$ be a packing in the plane with circular disks of radius at
most $1$.
Consider the set of disks $\D'$ in which each disk $D\in\D$ is enlarged by a
factor of $\Lambda = 1.00001$ from its center.
Then $\D'$ covers no square with side length $4$.
\end{thm}
There seems to be an overlooked case in the proof of this
theorem given in~\cite{kkmv-atpe-99}.
(Using the terminology of~\cite{kkmv-atpe-99}, one of the conditions that needs to be checked
in order to ensure that
$R_{n+1}$ is contained in~$R_n$ has been forgotten, see for
example the disk $D$ in Figure~11 of~\cite{kkmv-atpe-99}.)
We think that
the proof can be fixed, and moreover, the result can be proved for values somewhat
larger
than~1.00001, as the authors did not try to optimize the constant~$\Lambda$.
However, the case distinctions are very delicate, and we have not fully worked out
the details yet.
For these reasons, we state our main result depending on the value~$\Lambda$.

%
\begin{thm} \label{ImprovedLowerBound-thm}
Suppose Theorem~\ref{DiskPacking-thm} is true for a factor $\Lambda$ with
$\Lambda \leq 1.03$.
Then, the minimum geometric dilation $\Delta$ necessary to embed any finite
set of points in the plane satisfies
\[
  \Delta
  \geq
  \DeltaBound(\Lambda)
  :=
  \left(
    1
    +
    \left(\frac{\Lambda-1}{3}\right)^2
  \right)
  \frac{\pi}{2}
  .
\]
If Theorem~\ref{DiskPacking-thm} holds with $\Lambda = 1.00001$,
this results in
$\Delta \geq (1+(10^{-10})/9)\pi/2 > (1+10^{-11})\pi/2$.
\end{thm}
\begin{pf}
We first give an overview of the proof, and present the details afterwards.
Consider the set
$P := \left\{\, (x,y) \mid
x,y \in \left\{-9,-8,\ldots,9\right\}\;\right\}$
of 
grid points with integer coordinates in the square
$Q_1 := [-9,9]^2 \subset \realset^2$
, see Figure~\ref{PointSetAndConcatenated-fig}.
We use a proof by contradiction and assume that there exists a planar
connected graph $G$ that contains $P$ (as vertices or on its edges) and satisfies
$\delta(G) < \DeltaBound(\Lambda)$.
\begin{figure}[htb]%
  \begin{center}%
    \includegraphics{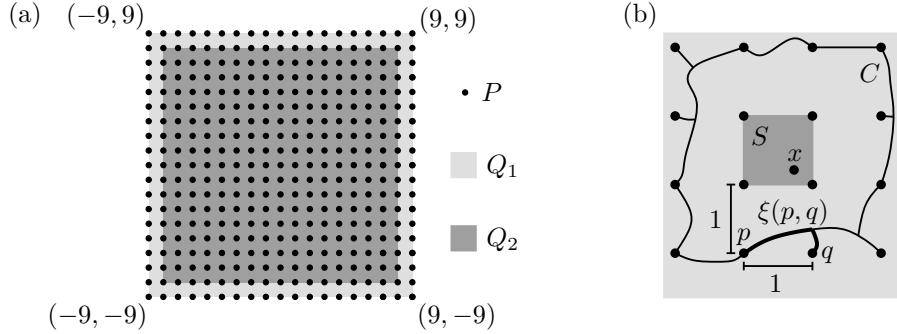}%
    \caption{\label{PointSetAndConcatenated-fig}
      (a)~The point set~$P$ and the square~$Q_2$
      which has to be covered by an embedding graph of low dilation.
      (b)~Every point $x \in Q_2$ is encircled by a cycle of length $\leq 12\delta(G)$.}%
  \end{center}%
\end{figure}
The idea of the proof is to show that if~$G$ attains such a low dilation,
it contains a collection~$\M$ of cycles with disjoint interiors
which cover the smaller square~$Q_2 := [-8,8]^2$.
We will choose $\M$ in such a way that the dilation of every cycle~$C \in \M$
will be bounded by
$
  \delta(C)
  \leq
  \delta(G)
  \leq
  \DeltaBound(\Lambda)
$.
The function~$\DeltaBound(\Lambda)$ is defined so that
we can easily apply Lemma~\ref{Ring-lem} to this situation.
In particular, $\Lambda \leq 1.03$ guarantees $\DeltaBound(\Lambda)\leq 1.0001 \pi/2$.
We derive that every $C$ has to be contained in
a $\Lambda$-ring.
As by Claim~\ref{BoundedLength-clm} the length of each cycle~$C \in \M$ will be bounded by~$8\pi$,
the inner disks of these rings have a radius~$r\leq 4$.
This is a contradiction to Theorem~\ref{DiskPacking-thm} (situation scaled by 4)
because the inner disks of the rings
are disjoint and their $\Lambda$-enlargements cover~$Q_2$.

We would like to use the cycles bounding the faces of~$G$
for~$\M$.
Indeed, $\delta(G)<\DeltaBound(\Lambda)<2$ implies that they cover~$Q_2$,
see Claim~\ref{BoundedLength-clm} below.
However, their dilation
could be bigger
than the dilation~$\delta(G)$ of the graph,
see Figure~\ref{NoShortcuts-fig}b.
There could be shortcuts in the exterior of~$C$,
i.e., the shortest path
between~$p,q \in C$ does not necessarily use~$C$.

Therefore, we have to find a different class of disjoint cycles
covering~$Q_2$ which do not allow shortcuts.
The idea is to consider for every point~$x$ in~$Q_2$
the shortest cycle of~$G$
such that~$x$ is contained in the open region bounded by the cycle.
These cycles are {\em non-crossing}, that is,
their enclosed open regions are either disjoint or one region contains the other
(Claim~\ref{NoCrossing-clm}).
If we define $\M$ to contain only the cycles of $\C$ which are maximal with respect to
inclusion of their regions, it provides all the properties we need.
We now present the proof in detail.
\begin{claim}\label{BoundedLength-clm}
Every point $x\in Q_2$ is enclosed by a cycle~$C$ of~$G$
of length at most~$8\pi$.
\end{claim}
\begin{pf}
We consider an arbitrary point $x \in Q_2$
and a grid square $S$ which contains~$x$
like shown in Figure~\ref{PointSetAndConcatenated-fig}b.
For every pair
$p,q$ of neighboring grid points
of~$P$,
let~$\xi(p,q)$ be
a shortest path in~$G$ connecting~$p$ and~$q$.
The length of each such path~$\xi(p,q)$ is bounded by $|\xi(p,q)| \le \delta_G(p,q)\cdot 1 < 2$.
Consider the closed curve $C$ obtained by concatenating the $12$
shortest paths between adjacent grid points on the boundary
of the $3 \times 3$ square around~$S$.
None of the shortest paths can enter $S$
because this would require a length bigger than 2.
Therefore~$C$ encloses but does not enter~$S$,
thus it also encloses~$x \in S$.
The total length of~$C$ is bounded by
$
\abs{C}
\leq 12\cdot\delta(G)
\leq 12 \DeltaBound(\Lambda)
\leq 12 \cdot 1.0001 (\pi/2)
< 8\pi
$.
\hfill$\qed$
\end{pf}
%
For any point $x \in Q_2$,
let $C(x)$ denote a \emph{shortest cycle} in $G$
such that $x$ is contained in the open region bounded by the cycle.
If the shortest cycle is not unique,
we pick one which encloses the smallest area.
It follows from Claim~\ref{NoCrossing-clm} below
that this defines the shortest cycle~$C(x)$ uniquely,
but this fact is not essential for the proof.
Obviously, $C(x)$ is a simple cycle (i.e., without
self-intersections).
Let $R(x)$ denote the open region enclosed by $C(x)$.
\begin{claim}\label{NoShortcuts-clm} For every $x \in Q_2$ we have:
    \begin{itemize}
      \item[{\rm(i)}] No shortest path of $G$ can cross $R(x)$.
      \item[{\rm(ii)}] Between two points $p,q$ on $C(x)$, there is always
        a shortest path on $C(x)$ itself.
    \end{itemize}
\end{claim}
\begin{figure}[htb]%
  \begin{center}%
    \includegraphics{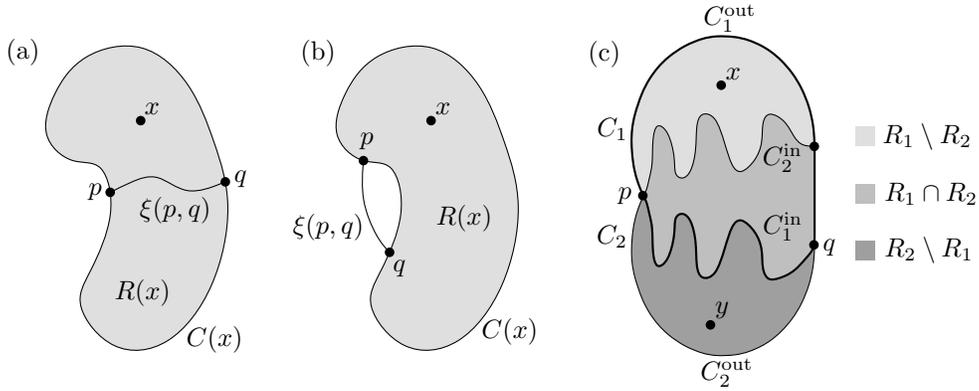}%
    \caption{\label{NoShortcuts-fig}
      Impossible situations: (a) a shortest path $\xi(p,q)$ crossing $R(x)$
      or (b) being a shortcut; (c) two crossing shortest cycles.%
    }%
  \end{center}%
\end{figure}%
\begin{pf}
  (i) Since every subpath of a shortest path is a shortest path, it suffices
  to consider a path~$\xi(p,q)$ between two points~$p,q$ on~$C(x)$ whose
  interior is completely contained in~$R(x)$, see Figure~\ref{NoShortcuts-fig}a.
  This path could replace one of the two arcs of~$C(x)$ between~$p$ and~$q$ and yield
  a better cycle enclosing~$x$, contradicting the definition of~$C(x)$.

  (ii)
  We have already excluded shortest paths which intersect~$R(x)$.
  We now exclude a path~$\xi(p,q)$ between two points~$p,q$ on~$C(x)$
  which runs outside~$C(x)$ and is strictly shorter than each of the two arcs
  of~$C(x)$ between
  $p$ and~$q$, see Figure~\ref{NoShortcuts-fig}b.
  If such a path existed, again, we could replace one of the
  arcs from~$C$ connecting~$p$ and~$q$ by $\xi(p,q)$
  and receive a shorter cycle encircling~$x$, a contradiction.
\hfill$\qed$
\end{pf}
As an immediate consequence of statement~(ii), we get:
\begin{claim}\label{BoundedDilation-clm}
  The dilation of every cycle~$C(x)$ is at most the dilation~$\delta(G)$ of
  the whole graph~$G$.
\end{claim}
This allows us to apply Lemma~\ref{Ring-lem} to the cycles~$C(x)$.
However, to obtain a packing, we still have to select a subset
of cycles with disjoint regions.
To this end we prove the following claim.
\begin{claim} \label{NoCrossing-clm}
For arbitrary points~$x,y \in Q_2$,
the cycles $C(x)$ and $C(y)$ are non-crossing,
i.e., $R(x)\cap R(y)=\emptyset \;\lor\; R(x)\subseteq R(y)
\;\lor\; R(x)\subseteq R(y)$.
\end{claim}
\begin{pf}
We argue by contradiction, see Figure \ref{NoShortcuts-fig}c.
Assume that the regions~$R_1$ and~$R_2$ of the shortest cycles~$C_1:=C(x)$
and~$C_2 := C(y)$ overlap, but none is fully contained inside the other.
This implies that their union~$R_1 \cup R_2$ is a bounded open connected set.
Its boundary~$\partial(R_1 \cup R_2)$ contains a simple cycle~$C$
enclosing~$R_1 \cup R_2$.

By the assumptions we know that there is a part~$C_1^{\textrm{in}}$
of~$C_1$ which connects two points~$p,q \in C_2$ and is, apart from its endpoints,
completely contained in~$R_2$.
Let~$C_1^{\textrm{out}}$ denote the other path on~$C_1$ connecting~$p$
and $q$.
By Claim~\ref{NoShortcuts-clm}(i),
at least one of the paths~$C_1^{\textrm{in}}$ or~$C_1^{\textrm{out}}$
must be a shortest path.
By Claim~\ref{NoShortcuts-clm}(ii),
$C_1^{\textrm{in}}$ cannot be a shortest path, since it intersects~$R_2$.
Hence, only~$C_1^{\textrm{out}}$ is a shortest path,
implying~$\abs{C_1^{\textrm{out}}} < \abs{C_1}/2$. 
Analogously, we can split $C_2$ into two paths~$C_2^{\textrm{in}}$ and~$C_2^{\textrm{out}}$
such that~$C_2^{\textrm{in}}$ is contained in~$R_1$, apart from its endpoints,
and~$\abs{C_2^{\textrm{out}}} < \abs{C_2}/2$. 

The boundary cycle~$C$ consists of parts of~$C_1$ and parts of~$C_2$.
It cannot contain any part of~$C_1^{\textrm{in}}$ or~$C_2^{\textrm{in}}$
because it intersects neither with~$R_1$ nor with~$R_2$.
Hence~$\abs{C} \leq \abs{C_1^{\textrm{out}}} + \abs{C_2^{\textrm{out}}}
< (\abs{C_1}+\abs{C_2})/2\le \max\{\abs{C_1},\abs{C_2}\}$.
Since~$C$ encloses~$x \in R_1$ and~$y \in R_2$, this contradicts
the choice of~$C_1 = C(x)$ or~$C_2 = C(y)$.
\hfill$\qed$
\end{pf}
Let~$\C$ be the set of shortest
cycles~$\C=\{\,C(x) \mid x \in Q_2 \,\}$,
and let~$\M \subset \C$ be the set
of \emph{maximal shortest cycles} with
respect to inclusion of their regions.
Claim~\ref{NoCrossing-clm} implies that these cycles have disjoint
interiors and that they cover~$Q_2$.
Claim~\ref{BoundedLength-clm} proves that their in-radius is bounded by $r \leq 4$.
By Claim~\ref{BoundedDilation-clm}, the dilation of
every cycle~$C \in \M$
satisfies~$\delta(C) \leq \delta(G) \leq \DeltaBound(\Lambda)$.
Like described in the beginning of this proof
we get a contradiction to Theorem~\ref{DiskPacking-thm}.
This completes the proof of the new lower bound
(Theorem~\ref{ImprovedLowerBound-thm}).
\hfill$\qed$
\end{pf}

\section{Closed Curves of Constant Halving Distance}
\label{ConstantHalvingDistance-sec}%
Closed curves of constant halving distance turn up naturally
if one wants to construct graphs of low dilation.
Lemma~\ref{DilationConvex-lem} shows that the dilation of a convex curve
of constant halving distance is attained by all its halving pairs.
Hence, it is difficult to improve (decrease) the dilation of such cycles,
because local changes decrease~$h$ or increase~$|C|$.

\begin{figure}[htbp]%
  \begin{center}%
      \begin{center}%
        \includegraphics{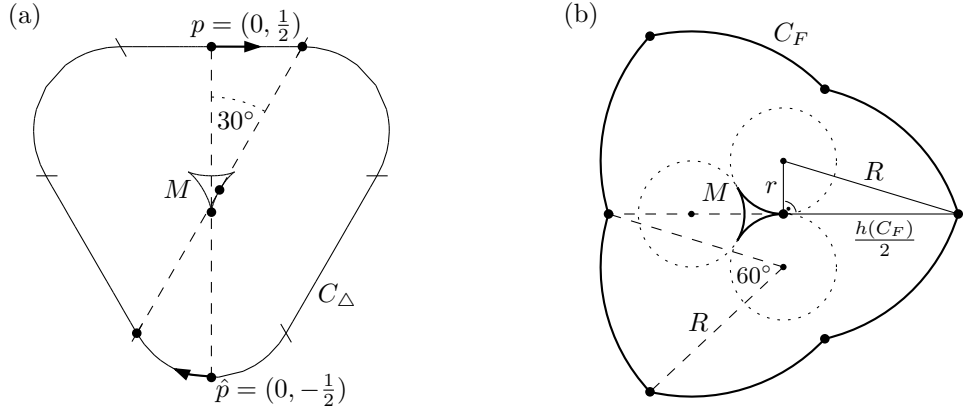}%
      \end{center}%
      \caption{\label{RoundedTriangleAndFlower-fig}
        The ``Rounded Triangle''~$C_\triangle$
        and the ``Flower'' $C_F$ 
        are curves of constant halving distance.}%
  \end{center}%
\end{figure}

This is the motivation for using
the curve of constant halving distance $C_F$ of
Figure~\ref{RoundedTriangleAndFlower-fig}b to construct
the grid of low dilation in Figure~\ref{ImprovedHexGrid-fig}.
Although the non-convex parts increase the length and thereby
the dilation of the small flowers $C_F$, they are useful in decreasing
the dilation of the big faces of the graph.

The proof of Lemma~\ref{Ring-lem}
shows that only curves with constant midpoint~$m(t)$
and constant halving distance can attain the
global dilation minimum of $\pi/2$.
Furthermore, it shows that only circles satisfy both conditions;
see also~\cite[Corollary 23]{egk-gdcpc-04ii},
\cite[Corollary 3.3]{acfgh-cmmke-02},
\cite{ks-dtk-98}, \cite{glp-smvr-81}.
What happens if only one of the conditions is satisfied?
Clearly, $m(t)$ is constant if and only if $C$ is centrally symmetric.
On the other hand, the class of closed curves of constant halving distance
is not as easy to describe.
One could guess --- incorrectly --- that it consists only of circles.
The ``Rounded Triangle''~$C_\triangle$ shown
in Figure~\ref{RoundedTriangleAndFlower-fig}a
is a convex counterexample, and
could be seen as an analog of the Reuleaux triangle,
the best-known representative of curves of constant width~\cite{cg-cbcw-83}.
It seems to be a somehow prominent example, because
two groups of the authors of this paper discovered it independently,
before finding it in a paper of H.~Auerbach~\cite[p.~141]{a-spmuc-38}
from 1938.

We construct $C_\triangle$ by starting with a pair of points
$p := (0,0.5)$ and $\hat{p} := (0,-0.5)$.
Next, we move $p$ to the right along a horizontal line.
Simultaneously, $\hat{p}$ moves to the left such that
the distance~$|p\hat{p}|=1$
is preserved and both points move with equal speed.
It can be shown that these conditions lead to a
differential equation whose solution
defines the path of $\hat{p}$ uniquely.
We move $p$ and $\hat{p}$ like this until the connecting
line segment $\edge{p\hat{p}}$
forms an angle of $30^\circ$ with the $y$-axis.
Next, we swap the roles of $p$ and $\hat{p}$.
Now, $\hat{p}$ moves along a line with the direction of its last movement,
and $p$ moves with equal speed on the unique curve which guarantees
$\abs{p\hat{p}}= 1$, until $\edge{p\hat{p}}$
has rotated with another $60^\circ$.
Again, we swap the roles of $p$ and $\hat{p}$ for
the next $60^\circ$ and so forth.
In this way we end up with six pieces of
equal length (three straight line segments and three curved
pieces) to build the Rounded Triangle~$C_{\triangle}$ depicted
in Figure~\ref{RoundedTriangleAndFlower-fig}a.
Note that the rounded pieces are not parts of circles.

The details of the differential equation and its solution are given in
Appendix~\ref{ParameterizationRoundedTriangle-sec}.
Here, we mention only that the perimeter of $C_\triangle$
equals $3 \ln 3$.
By Lemma~\ref{DilationConvex-lem}
this results in
\[
  \delta(C_{\triangle})
  =
  \pseudofrac{|C_{\triangle}|}{(2h(C_{\triangle}))}
  =
  \frac{3}{2} \ln 3
  \approx
  1.6479
  \;.
\]
The midpoint curve of $C_\triangle$ is built from six congruent pieces that
are arcs of a tractrix,
as we will discuss in the end of this section.
First, we give a necessary and sufficient condition for
curves of constant halving distance (not necessarily convex).
\begin{thm}\label{ConstantHalvingDistance-thm}
  Let $C$ be a planar closed curve,
  and let $c: [0,\abs{C}) \rightarrow C$
  be an arc-length parameterization.
  Then, the following two statements are equivalent:
  \begin{enumerate}
    \item\label{CurveConditions}
      If $c$ is differentiable in $t$ and in $t+\halfCtext$,
      $\dot{m}(t) \neq 0$, and $\dot{c}^*(t) \neq 0$, then
      the halving chord $c(t)c(t+\halfCtext)$ is tangent to the midpoint curve at $m(t)$.
      And if the midpoint stays at $m \in \realset^2$
      on a whole interval $(t_1,t_2)$,
      the halving pairs are located on the circle
      with radius $\pseudofrac{h(C)}{2}$
      and center point $m$.
    \item\label{ConstantHalvingDistance}
      The closed curve $C$ is a cycle of constant halving distance.
  \end{enumerate}
\end{thm}
\begin{pf}
  ``{\em \ref{ConstantHalvingDistance}. $\Rightarrow$ \ref{CurveConditions}.}''

  Let $C$ have constant halving distance.
  If $c$ is differentiable in $t$ and $t+\halfCtext$,
  $c^*$ and $m$ are differentiable in $t$.
  And due to $\abs{c^*}\equiv\pseudofrac{h(C)}{2}$
  it follows %
  that $\dot{c}^*(t)$ must be orthogonal to $c^*(t)$ which can be shown by
  \begin{equation}\label{CStarDerivativeOrthogonal-equ}
    0
    =
    \frac{\intd}{\intd t}
    |c^*(t)|^2
    =
    \frac{\intd}{\intd t}
    \left<c^*(t),c^*(t)\right>
    =
    2\left<c^*(t),\dot{c}^*(t)\right>
    \;.
  \end{equation}
  On the other hand,
  we have already seen in \thetag{\ref{MCStarOrthogonal-equ}}
  of Section~\ref{MAndCStar-sec} that $\dot{m}(t)$
  and $\dot{c}^*(t)$ are orthogonal.
  Hence, $\dot{m}(t) \neq 0 \neq \dot{c}^*(t)$ implies $\dot{m}(t) \parallel c^*(t)$
  and the first condition of \ref{CurveConditions}. is proven.
  The second condition follows trivially from $c(t)=m(t)+c^*(t)$.

  \vspace{1 em}
  \pagebreak[1]
  \noindent``{\em \ref{CurveConditions}. $\Rightarrow$ \ref{ConstantHalvingDistance}.}''

  Let us assume that~\ref{CurveConditions}. holds.
  We have to show that~$\abs{c^*(t)}$ is constant.

  First, we consider an interval~$(t_1,t_2)\subseteq [0,|C|)$,
  where~$m(t)$ is constant~($=m$)
  and the halving pairs are located on
  a circle with radius~$\pseudofrac{h(C)}{2}$ and center~$m$.
  This immediately implies that~$\abs{c^*}$ is constant on $(t_1,t_2)$.

  If~$(t_1,t_2)\subseteq [0,|C|)$ denotes an interval
  where~$|c^*(t)|=0$, then obviously~$\abs{c^*}$ is constant.

  Now, let~$(t_1,t_2)\subseteq [0,|C|)$ be an open interval
  where~$c(t)$
  and~$c\left(t+\frac{\abs{C}}{2}\right)$
  are differentiable and~$\dot{m}(t)\neq0$ and~$\dot{c}^*(t) \neq 0$
  for every~$t \in (t_1,t_2)$.
  We follow the proof
  of ``{\em \ref{ConstantHalvingDistance}. $\Rightarrow$ \ref{CurveConditions}.}''
  in the opposite direction.
  Equation \thetag{\ref{MCStarOrthogonal-equ}} shows that
  $\dot{c}^*(t) \perp \dot{m}(t)$ and the first condition of 1.
  gives $c^*(t) \parallel \dot{m}(t)$.
  Combining both statements results in $\dot{c}^*(t) \perp c^*(t)$
  which by~\thetag{\ref{CStarDerivativeOrthogonal-equ}}
  yields that $|c^*(t)|$ is constant.

  The range~$\left[0,\halfCtext\right)$ can be divided into countably many
  disjoint intervals $[t_i,t_{i+1})$ where~$m$ and~$c^*$ are differentiable
  on the open interval~$(t_i,t_{i+1})$,
  and one of the three conditions
  $\dot{m}(t)=0$, $\dot{c}^*(t)=0$ or $\dot{m}(t) \neq 0 \neq \dot{c}^*(t)$
  holds for the whole interval~$(t_i,t_{i+1})$.
  We have shown that~$\abs{c^*}$ must be constant on all these open intervals.
  Thus, due to~$c^*$ being continuous on~$[0,\halfCtext)$,
  $\abs{c^*}$ must be globally constant.
\hfill$\qed$
\end{pf}
The theorem shows that curves of constant halving distance can
consist of three types of parts; parts where the halving chords
lie tangentially to the midpoint curve, circular arcs
of radius $h(C)/2$, and parts where $\dot c^*(t)=0$ and the halving pairs are only moved
by the translation due to $m$,
i.e., for every $\tau_1$ and $\tau_2$ within such a part
we have $c(\tau_2)-c(\tau_1) = c(\tau_2+\halfCtext) - c(\tau_1+\halfCtext)
= m(\tau_2)-m(\tau_1)$.
For convex cycles of constant halving pair distance, the translation parts
cannot occur:
\begin{lem}\label{NoTranslationParts-lem}
  Let $C$ be a closed convex curve of constant halving distance.
  Then there exists no non-empty interval $(t_1,t_2) \subset [0,|C|)$
  such that $c^*$ is constant on $(t_1,t_2)$.
\end{lem}
\begin{pf}
Assume that $c^*$ is constant on $(t_1,t_2)$ and
choose $s_1,s_2$ with $t_1<s_1<s_2<t_2$ and $s_2<s_1+|C|/2$.
If the four points
$p_1=c(s_1)$,
$p_2=c(s_2)$,
$p_3=c(s_2+|C|/2)$,
$p_4=c(s_1+|C|/2)$
do not lie on a line, they form a parallelogram in which $p_1p_4$ and
$p_2p_3$ are parallel sides.
However, these points appear on $C$ in the cyclic order $p_1p_2p_4p_3$, which
is different from their convex hull order $p_1p_2p_3p_4$ (or its reverse),
a contradiction.
The case when the four points lie on a line $\ell$ can be dismissed easily
(convexity of $C$ implies that the whole curve $C$
would have to lie on $\ell$, but then $C$ could not be a curve of
constant halving distance).
\hfill$\qed$
\end{pf}

\begin{figure}[htbp]
        \begin{center}%
          \includegraphics{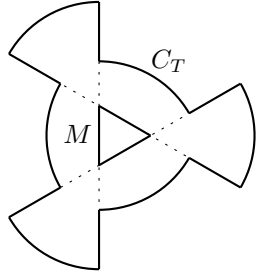}%
          \caption{\label{CycleWithTranslation-fig} $C_T$
            consists of translation parts and circular arcs.}%
        \end{center}%
\end{figure}

If we drop the convexity condition, there are easy examples of cycles
of constant halving distance with translation parts, see for example
Figure~\ref{CycleWithTranslation-fig}.
(Of course, the translation parts need not be line segments.)
However, the proof extends to
the case when the halving chords go through the interior of~$C$. (When
halving chords may touch the boundary of~$C$, then the case where
$p_1p_2p_4p_3$ lie on a line persists, as in
Figure~\ref{CycleWithTranslation-fig}.)

With the convexity requirement, we have proved that the halving chord is
always tangent to the midpoint curve in all points where the midpoint curve
has a tangent: the midpoint curve is the \emph{envelope} of the
halving segments. This is the viewpoint from which 
Zindler  started his investigations
in 1921: he was interested in envelopes of all sorts of halving chords
(halving the perimeter, the area, etc., and not necessarily of constant
length) and other classes of chords with special properties.
 For curves of constant halving distance, he observed the following
interesting fact (without stating it explicitly as a
theorem):
\begin{thm} \textrm{\rm (Zindler~\cite[Section~7]{z-ukg-21})}
For a convex closed curve in the plane,
the following statements are equivalent.
\begin{enumerate}
\item All halving chords have the same length.
\item All chords halving the \emph{area} have the same length.
\item Halving chords and area-halving chords coincide.
\end{enumerate}
\end{thm}
Following Auerbach~\cite{a-spmuc-38}, these curves are consequently
called \emph{Zindler curves}.  (The theorem holds also for non-convex
curves as long as all halving chords go through the interior of~$C$,
and it remains even true for chords that divide the area or perimeter
in some arbitrary constant ratio.)  Auerbach~\cite{a-spmuc-38} noted
than these curves are related to a problem of S.~Ulam about floating
bodies.  Ulam's original problem~\cite[Problem~19,
pp.~90--92]{m-sbmsc-82}, \cite[Section~A6, p.~19--20]{cfg-upg-91},
\cite[pp.~153--154]{s-ms-69}, which is still unsolved in its generality,
asks if there is a homogeneous body different from a spherical ball that can
float in equilibrium in any orientation.
A two-dimensional version can be formulated as follows: consider a homogeneous
cylindrical log of density $\frac12$ whose cross-section is a convex
curve~$C$, floating in water (which has density~1).
This log will float in
equilibrium \emph{with respect to rolling} in every \emph{horizontal}
position if and only if $C$ is a Zindler curve.

Moreover, Auerbach
 relates these curves to curves of
constant width~\cite[\S~7]{cg-cbcw-83}:
\begin{thm}\label{th-auerb}\textrm{\rm (Auerbach~\cite{a-spmuc-38})}
If a square $PQRS$ moves rigidly in the plane such the diagonal $PR$
traces all halving chords of a Zindler curve $C$
the endpoints of the other diagonal $QS$ trace out a
curve $D$ of constant width, enclosing the same area as~$C$.
\end{thm}
Conversely, one can start from any curve of constant width.
Taking point pairs $Q$ and $S$ with parallel tangents as diagonals of squares
$PQRS$, the other diagonal $PR$ will generate a curve of
constant halving distance, not necessarily convex.
Theorem~\ref{th-auerb} remains true for a rhombus $PQRS$ as long as
$QS$ is it not too short, relatively to $PR$. (Of course, the area of $D$ will then be different.)

\begin{figure}[htbp]%
        \begin{center}
          \includegraphics{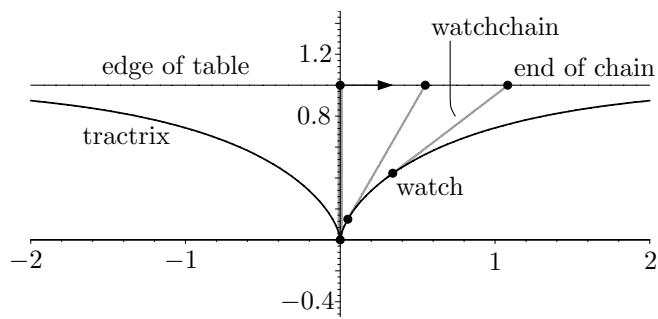}%
          \caption{\label{Tractrix-fig} The tractrix, the curve of a
            watch on a table towed with its watchchain
            (the curve is symmetric about the $y$-axis).}%
        \end{center}
\end{figure}

Now we show that the midpoint curve of the
Rounded Triangle~$C_\triangle$ is built from six tractrix pieces.
The tractrix is illustrated in Figure~\ref{Tractrix-fig}.
A watch is placed on a table, say at the origin $(0,0)$
and the end of its watchchain of length $1$ is pulled along the horizontal edge
of the table starting at $(0,1)$, either to the left or to the right.
As the watch is towed in the direction of the chain, the chain is
always tangent to the path of the watch, the tractrix.

From the definition it is clear that the midpoint curve of
the cycle $C_\triangle$
consists of such tractrix pieces, scaled by $1/2$,
because by definition and Theorem~\ref{ConstantHalvingDistance-thm}
its halving chords are always tangent to the midpoint curve,
always one of the points of these pairs is moving on a straight line,
and its distance to the midpoint curve stays~$1/2$.
A parameterization not depending on the midpoint curve is
analyzed in Appendix~\ref{ParameterizationRoundedTriangle-sec}.

As mentioned in the introduction, Zindler~\cite[Section~7.b]{z-ukg-21} discovered some
convex curves of constant halving distance already in 1921.
He restricted his analysis to curves whose midpoint curves have three cusps
like the one of the rounded triangle.
His example uses as the midpoint curve
the path of a point on a circle of radius $1$ rolling inside of
a bigger circle of radius $3$ (a hypocycloid). 
To guarantee convexity, one then needs $h \geq 48$.

\section{Relating Halving Distance to other Geometric Quantities}%
\label{hAndOther-sec}%
One important topic in convex geometry is the
relation between different geometric quantities of convex bodies
like area~$A$ and diameter~$D$.
Scott and Awyong~\cite{sa-ics-00} give a short survey
of basic inequalities in $\realset^2$.
For example, it is known that $4A \leq \pi D^2$,
and equality is attained only by circles,
the so-called {\em extremal set} of this inequality.

In this context the minimum and maximum halving distance
$h$ and~$H$ give rise to some new interesting questions,
namely the relation to other basic quantities like the width $w$.
As the inequality $h \leq w$ is immediate from definition,
the known upper bounds on~$w$ hold for~$h$ as well.
However, although the original inequalities
relating~$w$ to perimeter, diameter, area, inradius and circumradius are tight
(see~\cite{sa-ics-00}), not all of them are also tight for~$h$.
One counterexample ($A \geq w^2/\sqrt{3} \geq h^2/\sqrt{3}$)
will be discussed in the following subsection.

\comment{%
Halving chords have been investigated for a long time.
Already in 1921,
Zindler~\cite{z-ukg-21} deals with many different kinds of chords, 
such as our (perimeter) halving chords and area halving chords.
Among other, he shows the existence of convex curves that
are not centrally symmetric, and for which perimeter halving chords
and area halving chords are the same in each direction.
}

\subsection{Minimum Halving Distance and Area}
\label{hAndA-sec}
Here, we consider the relation of the minimum halving distance~$h$
and the area~$A$.
Clearly, the area can get arbitrarily big while $h$ stays constant.
For instance this is the case for a rectangle of smaller side length~$h$ where
the bigger side length tends to infinity.

How small the area~$A$ can get for closed convex curves
of minimum halving distance~$h$?
A first answer~$A \geq h^2/\sqrt{3}$
is easy to prove,
because it is known~\cite[ex. 6.4, p.221]{yb-cf-61}
that~$A \geq w^2/\sqrt{3}$,
and we combine this with $w \geq h$.
Still, $1/\sqrt{3} \approx 0.577$ is not
the best possible lower bound to $A/h^2$,
since the equilateral triangle is the only
closed curve attaining~$A = w^2/\sqrt{3}$
and its width~$w=\pseudofrac{\sqrt{3}}{2}\approx 0.866$ (for side length 1)
is strictly bigger than its minimum halving distance~$h=3/4=0.75$.
We do not know the smallest possible value of~$A / h^2$.
\begin{figure}[htbp]%
  \begin{center}%
      \begin{center}%
        \includegraphics{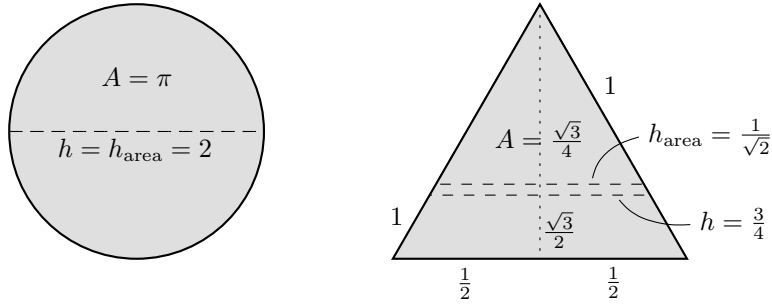}%
      \end{center}%
      \caption{\label{EquilateralTriangle-fig} The equilateral triangle
        has a smaller ratio
        $A/h^2=\pseudofrac{4}{(3\sqrt{3})}\approx 0.770$ and
        a bigger ratio
        $A/h_{\mathrm{area}}^2=\pseudofrac{\sqrt{3}}{2}\approx 0.866$
        than the circle ($\pi/4\approx 0.785$).}
  \end{center}%
\end{figure}

We can also consider chords bisecting the
area of planar convex sets instead of halving their perimeter,
see~\cite[Section~A26, p.~37]{cfg-upg-91} about ``dividing up a piece
of land by a short fence''.
Let $h_{\mathrm{area}}$ be the minimum area-halving distance.
Santal{\'o} asked%
\footnote{We would like to thank Salvador Segura Gomis
for pointing this out.}
if
$A \geq (\pseudofrac{\pi}{4}) h_{\mathrm{area}}^2$.
 Equality is attained by a circle.

\begin{figure}[htbp]%
  \begin{center}%
      \begin{minipage}[t]{5cm}%
        \begin{center}%
          \includegraphics{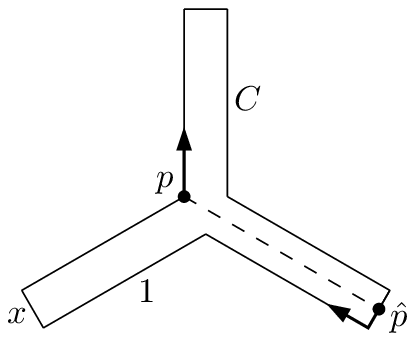}%
          \caption{\label{CycleOfSmallArea-fig} By decreasing $x$ we can make
            the area of this closed curve arbitrarily small while $h$ stays
            bounded.}%
        \end{center}%
      \end{minipage}%
      \hfill
      \begin{minipage}[t]{8.0cm}%
        \begin{center}%
          \includegraphics{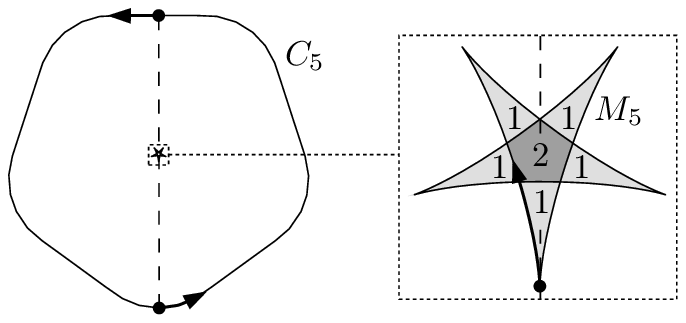}%
          \caption{\label{MRounded5gon-fig} The midpoint curve of
            a rounded pentagon, constructed analogously to $C_\triangle$
            of Figure~\ref{RoundedTriangleAndFlower-fig},
has 5 cusps. It
            contains regions with winding number 2
            and regions with winding number 1.}%
        \end{center}%
      \end{minipage}%
  \end{center}%
\end{figure}

Now going back to {\em perimeter} halving distance, does the circle attain
smallest area for a given~$h$? Already the equilateral
triangle gives a counterexample,
$A/h^2 = \frac{\sqrt{3}}{4}/\frac{9}{16}
\approx 0.770 < 0.785 \approx \frac{\pi}{4}$,
see Figure~\ref{EquilateralTriangle-fig}.
Still, we do not know if the equilateral triangle is the convex cycle
minimizing $A/h^2$.
On the other hand, clearly, $A/h^2$ can become arbitrarily small
if we drop the convexity condition, see Figure~\ref{CycleOfSmallArea-fig}.
Not only does the equilateral triangle attain a smaller ratio
$A/h^2$ than the circle.
In fact, the circle maximizes~$A/h^2$
among the curves of constant halving distance.
The following formula was proved by Zindler~\cite[Theorem~24]{z-ukg-21}
for the special case of midpoint curves
with three cusps.
Here, we give a proof for the general case. 
\begin{lem}\label{ConstantHDArea-lem}
  If~$C$ is a convex cycle of constant halving distance~$h$,
  its area satisfies~$A = (\pseudofrac{\pi}{4}) h^2 - 2A(M)$
  where $A(M)$ denotes the area bounded by the
  midpoint curve $M$.
  In $A(M)$, the area of any region encircled several times by $M$
  is counted with the multiplicity of the corresponding winding number,
  see Figure~\ref{MRounded5gon-fig} for an example.
  In particular, $A \leq (\pseudofrac{\pi}{4}) h^2$.
\end{lem}
\begin{pf}
  It suffices to analyze the case $h=2$.
  And we assume that $c$ is oriented counterclockwise.
  Then, we can consider the following parameterizations of $C^*$ and $M$
  in terms of the angle $\alpha$ of the current halving chord.
  \[
    c^*(\alpha)=\vect{\cos \alpha}{\sin \alpha}
    ,\qquad
    \dot m(\alpha)= v(\alpha) \cdot \vect{\cos \alpha}{\sin \alpha}
  \]
  Such parameterizations exist for some function $v:[0,2\pi)\to\realset$
  because of Theorem~\ref{ConstantHalvingDistance-thm}
  and Lemma~\ref{NoTranslationParts-lem}.
  Periodicity requires $m(\alpha)=m(\alpha+\pi)$
  and thereby $v(\alpha+\pi)=-v(\alpha)$.
  We say that $v$ is {\em antiperiodic} in $\pi$.
  If~$x(\alpha)$ denotes the $x$-coordinate of~$c(\alpha)=c^*(\alpha)+m(\alpha)$
  and $y(\alpha)$ is the $y$-coordinate, we get
  \begin{eqnarray*}
    A
    & = &
    \int x \;\intd y
    =
    \int_0^{2\pi} c_x(\alpha) \dot{c}_y(\alpha) \intd \alpha
    \\ & \stackrel{\scriptstyle c=c^*+m}{=} &
    \int_0^{2\pi}
    \left(
      \cos \alpha + \int_0^\alpha v(\beta)\cos \beta \,\intd \beta
    \right)
    \left(
      \cos \alpha + v(\alpha)\sin\alpha
    \right)
    \intd \alpha
    \\ & = &
    \int_0^{2\pi} \cos^2 \alpha \,\intd \alpha
    \;+\;
    \int_0^{2\pi}
      \underbrace{v(\alpha)}_{\textrm{\scriptsize antip.}}\;
      \underbrace{\cos \alpha}_{\textrm{\scriptsize antip.}}\;
      \underbrace{\sin \alpha}_{\textrm{\scriptsize antip.}}\;
      \intd \alpha
    \\ &&
    \;+\;
    \int_0^{2\pi}
      \underbrace{\cos \alpha}_{\textrm{\scriptsize antip.}}\;
      \underbrace{
        \left(
          \int_0^\alpha v(\beta) \cos \beta \,\intd \beta
        \right)
      }_{\mbox{\scriptsize $=m_x(\alpha)$, periodic in $\pi$}}\;
      \intd \alpha
    \\ &&
    \;+\;
    \int_0^{2\pi}
    v(\alpha) \sin \alpha
    \left(
      \int_0^\alpha v(\beta) \cos \beta \,\intd \beta
    \right)
    \intd \alpha
    \\ 
    & = &
    \pi
    +
    2
    \iint\limits_{0 \le \beta \le \alpha \le \pi}
    v(\alpha)v(\beta)\sin{\alpha}\cos{\beta}
    \,\intd\beta\;\intd\alpha
    .
  \end{eqnarray*}
  The two integrals with marked terms disappear in the above equations
  because their integrands are antiperiodic in $\pi$.
  For the area of $M$ we get analogously
  \[
    A(M)
    =
    - \iint\limits_{0 \le \beta \le \alpha \le \pi}
    v(\alpha)v(\beta)\sin{\alpha}\cos{\beta}
    \,\intd\beta \;\intd\alpha
    .
  \]
  The negative sign here comes from the fact that $M$ is
  traversed in the opposite orientation from $C$.
\hfill$\qed$
\end{pf}
We have shown that the circle is the convex cycle of constant halving distance attaining maximum area.
But which convex cycle of constant halving distance attains minimum area?
We conjecture (in accordance with
Auerbach~\cite[p.~138]{a-spmuc-38})
 that the answer is the Rounded Triangle~$C_\triangle$.
Lemma~\ref{ConstantHDArea-lem} helps us to calculate its area
$A(C_\triangle)$.
The tractrix construction of the midpoint curve $M$
makes it possible to get a closed form for $A(M)$.
The value of the resulting expression for the
area of~$C_\triangle$ is
\[
  A(C_\triangle) = \sqrt3\cdot\left[\log 3 - (\log^2 3)/8 -
    1/2\right]\cdot h^2
 \approx 0.7755\cdot h^2
  \;.
\]
This value (as well as the length of $C_\triangle$) has also been
obtained
by Auerbach~\cite{a-spmuc-38}, with a different method.
Auerbach also proposed $C_\triangle$ as a candidate for a Zindler
curve
of maximum length.

Zindler's curve~\cite[Section~7.b]{z-ukg-21}, which we described in the end of
Section~\ref{ConstantHalvingDistance-sec}, has 
area~$((\pi/4)-(1/24)^2\pi)h^2 \approx 0.780 h^2$.
Both constant factors are smaller than
$\pi/4\approx 0.7854\ldots$, thus providing a negative answer
to the above-mentioned question of Santal{\'o} whether
$A \geq (\pseudofrac{\pi}{4}) h_{\mathrm{area}}^2$.
It would be interesting to know whether $C_\triangle$ has the smallest
area among all convex curves with a given minimum area-halving distance
$h_{\mathrm{area}}$.

Auerbach~\cite{a-spmuc-38} constructed another, non-convex, curve of
constant halving width, based on the same tractrix construction as the
rounded triangle $C_\triangle$, but consisting of only two straight
edges and two smooth arcs, forming the shape of a heart.  Approximate
renditions of both curves can be found in several problem collections
\cite[Fig.~19.1, p.~91]{m-sbmsc-82}, \cite[Fig.~A9,
p.~20]{cfg-upg-91}, \cite[Figs.~179 and 180, pp.~153--154]{s-ms-69}
and other popular books, but they appear to use circular arcs instead
of the correct boundary curves, given by the parametrization
\thetag{\ref{eq:parame}} in
Appendix~\ref{ParameterizationRoundedTriangle-sec}.

\subsection{Minimum Halving Distance and Width}
In order to achieve a lower bound to~$h$ in terms of~$w$, we examine
the relation of both quantities to the area~$A$ and the diameter~$D$.
The following inequality was first proved by Kubota~\cite{k-uee-23} in 1923
and is listed in~\cite{sa-ics-00}.
\begin{thm}{\em (Kubota~\cite{k-uee-23})}
\label{InequalityAreaDiameterWidth-thm} If $C$ is a convex curve, then $ A \geq {Dw}/{2} $.
\end{thm}
We will combine this known inequality with the following new result.
\begin{thm} \label{InequalityAreaDiameterHD-thm}
  If $C$ is a convex curve, then $A \leq Dh$.
\end{thm}
\begin{pf}
Without loss of generality we assume that a halving
chord $p\hat{p}$ of minimum length~$h$ lies on the $y$-axis,
$p$ on top and $\hat{p}$ at the bottom,
see Figure~\ref{InequalityAhD-fig}.
Let $C_-$ be the part of $C$ with negative $x$-coordinate
and let $C_+ := C \setminus C_-$ be the remaining part.
We have $|C_-|=|C_+|=|C|/2$ because $p\hat{p}$ is a halving chord.
In Figure~\ref{InequalityAhD-fig}, obviously $|C_-| = |C_+|$ does not
hold, but this is only to illustrate our proof by contradiction.
\begin{figure}[htbp]%
  \begin{center}%
      \begin{center}%
        \includegraphics{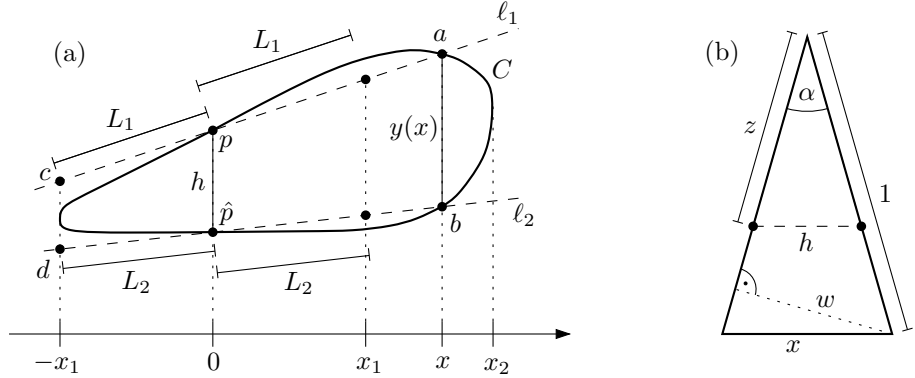}%
      \end{center}%
      \caption{\label{InequalityAhD-fig}
        (a)~Proving by contradiction that $y(x)\leq h$
          for every $x$ in $[x_1,x_2]$.
        (b)~In a thin isosceles triangle $h/w \searrow 1/2$ if $\alpha \to 0$.}%
  \end{center}%
\end{figure}
Let $-x_1$ and $x_2$ denote the minimum and maximum $x$-coordinate of $C$.
Note that $x_1$ has a positive value.
We assume that $x_2>x_1$.
Otherwise we could reflect $C$ at the $y$-axis.
Let $y(x)$ be the length of the vertical line segment of $x$-coordinate $x$
inside~$C$, for every $x \in [-x_1,x_2]$.
These definitions result in $x_1+x_2 \leq D$
and $A = \int_{-x_1}^{x_2} y(x)\,\intd x$.
Furthermore, if we take the convex hull
of the vertical segment with $x$-coordinate $x$ and
length $y(x)$ and the vertical segment with $x$-coordinate $-x$
and length $y(-x)$, then its intersection with the $y$-axis has
length $(y(-x)+y(x))/2$.
By convexity it must be contained in the line segment $p\hat{p}$
of length~$h$.
This implies
\begin{equation}
  \label{YMeanValueProperty-equ}
  \forall x \in [0,x_1]:\;\;
  y(-x)+y(x) \leq 2h
  \;.
\end{equation}
As a next step, we want to show that
\begin{equation}
  \label{YUpperBoundEnd-equ}
  \forall x \in [x_1,x_2]:\;\;
  y(x)\leq h
  \;.
\end{equation}
We assume that $y(x) > h$.
Let $ab$ be the vertical segment of $x$-coordinate $x$ inside $C$,
$a$ on top and $b$ at the bottom.
Then, we consider the lines $\ell_1$ through $p$ and $a$
and $\ell_2$ through $\hat{p}$ and $b$.
Let $L_1$ ($L_2$) be the length of the piece of $\ell_1$ ($\ell_2$)
with $x$-coordinates in~$[0,x_1]$.
By construction they are equal to the corresponding lengths
in the $x$-interval~$[-x_1,0]$.
Let $c$ and $d$ be the points with $x$-coordinate $-x_1$ on $\ell_1$, $\ell_2$ respectively.
Then, by the convexity of $C$, we have $|C_-| \leq L_1 + L_2 + |cd| \leq L_1 + L_2 + h
< L_1 + L_2 + y(x) \leq |C_+|$.
This contradicts to $p\hat{p}$ being a halving chord,
and the proof of (\ref{YUpperBoundEnd-equ}) is complete.

Now we can plug everything together and get
\begin{eqnarray*}
  A
  & = &
  \int\limits_{-x_1}^{x_2} y(x) \,\intd x
  =
  \int\limits_{0}^{x_1} \left(y(-x)+y(x)\right) \,\intd x
  \;
  +
  \;
  \int\limits_{x_1}^{x_2} y(x) \,\intd x
  \\&\stackrel{\textrm{\scriptsize (\ref{YMeanValueProperty-equ}),\;(\ref{YUpperBoundEnd-equ})}}{\leq}&
  x_1\cdot2h + (x_2-x_1)h
  =
  (x_1+x_2)h
  \leq Dh
  \;.
\end{eqnarray*}
\hfill$\qed$
\end{pf}
Finally, we achieve the desired inequality relating $h$ and~$w$.
\begin{thm}\label{InequalityHalvingDistanceWidth-thm}
If $C$ is a convex curve, then $h \geq w/2$.
This bound is tight.
\end{thm}
\begin{pf}
The inequality follows directly from
Theorem~\ref{InequalityAreaDiameterWidth-thm}
and Theorem~\ref{InequalityAreaDiameterHD-thm}.
To see that the bound is tight, consider a thin isosceles triangle
like in Figure~\ref{InequalityAhD-fig}b.
If $h$ is the minimum halving distance,
we have $2z = |C|/2 = 1+x/2 = 1+\sin(\alpha/2)$, thus
$h=2z\sin\frac{\alpha}{2}=(1+\sin(\alpha/2))\sin(\alpha/2)$.
On the other hand, the width is given by
$w = \sin \alpha = 2\sin(\alpha/2)\cos(\alpha/2)$, therefore
$h/w = (1+\sin(\alpha/2))/(2\cos(\alpha/2)) \searrow 1/2$
for $\alpha \to 0$.
\hfill$\qed$
\end{pf}
Theorem~\ref{InequalityHalvingDistanceWidth-thm} can be also proved
directly by using arguments analogous to the proof of
Theorem~\ref{InequalityAreaDiameterHD-thm}.
But we think that Theorem~\ref{InequalityAreaDiameterHD-thm}
is of independent interest.

\section{Dilation Bounds}
\label{DilationBounds-sec}%
\subsection{Upper Bound on Geometric Dilation}
\begin{figure}[htbp]%
  \begin{center}%
      \begin{center}%
        \includegraphics{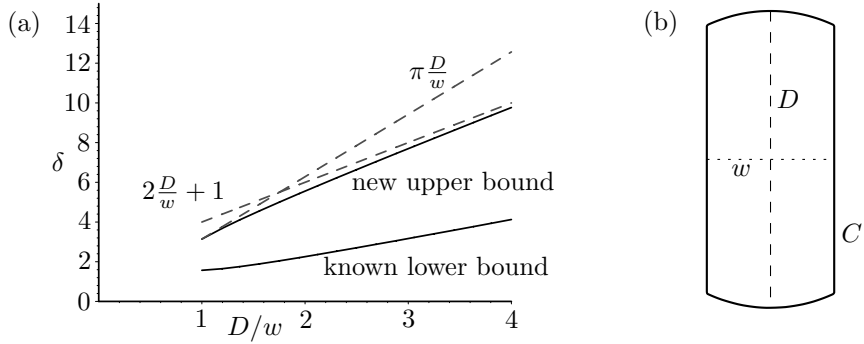}%
      \end{center}%
      \caption{\label{DilationBoundsCycle-fig}
        (a)~A plot of the new upper bound on geometric dilation
        and the known lower bound.
        (b)~The extremal set of Kubota's
        inequality~\thetag{\ref{LengthUpperBound-equ}}.}%
  \end{center}%
\end{figure}

Our Theorem~\ref{InequalityHalvingDistanceWidth-thm}
leads to a new upper bound depending only on the ratio $D/w$.
This 
complements the lower bound
\begin{equation}\label{LowerDilationBoundDw-equ}
  \delta(C) \geq \arcsin{\frac{w}{D}} + \sqrt{\left(\frac{D}{w}\right)^2-1}
\end{equation}
of Ebbers-Baumann et al.~\cite[Theorem 22]{egk-gdcpc-04ii}.
The new upper bound is stated in the following theorem
and plotted in Figure~\ref{DilationBoundsCycle-fig}a.
\begin{thm} \label{UpperDilationBound-lem}
If $C$ is a convex curve, then
\[
  \delta(C)
  \leq
  2
  \left(
    \frac{D}{w}
    \arcsin{\frac{w}{D}}
    +
    \sqrt{\left(\frac{D}{w}\right)^2-1}
  \right)
  .
\]
\end{thm}
\begin{pf}
  Kubota~\cite{k-uee-23} (see also~\cite{sa-ics-00}) showed that
  \begin{equation}
    \label{LengthUpperBound-equ}
    |C| \leq 2D\arcsin{\frac{w}{D}}+2\sqrt{D^2-w^2}
    \;.
  \end{equation}
  Combining this with Theorem~\ref{InequalityHalvingDistanceWidth-thm}
  and Lemma~\ref{DilationConvex-lem} yields
  \[
    \delta(C)
    \stackrel{\textrm{\scriptsize 
Lemma~\ref{DilationConvex-lem}}}{=}
    \frac{|C|}{2h}
    \stackrel{\textrm{\scriptsize Theorem~\ref{InequalityHalvingDistanceWidth-thm}}}{\leq}
    \frac{|C|}{w}
    \stackrel{\textrm{\scriptsize (\ref{LengthUpperBound-equ})}}{\leq}
    2
    \left(
      \frac{D}{w}
      \arcsin{\frac{w}{D}}
      +
      \sqrt{\left(\frac{D}{w}\right)^2-1}
    \right)
    \;.
  \]
\hfill$\qed$
\end{pf}
The isoperimetric inequality $|C| \leq D \pi$ and
the inequality $|C| \leq 2D + 2w$ lead to the
slightly bigger but simpler dilation bounds
$\delta(C) \leq \pi \frac{D}{w}$
and $\delta(C) \leq 2\left(\frac{D}{w}+1\right)$,
see Figure~\ref{DilationBoundsCycle-fig}a.
But even the dilation bound of Theorem~\ref{UpperDilationBound-lem} is not tight
because (\ref{LengthUpperBound-equ}) becomes an equality only for
curves which result from the intersection of a circular disk of
diameter $D$ with a parallel strip of width $w$,
see Figure~\ref{DilationBoundsCycle-fig}b.
For these curves we have $h=w$ due to their central symmetry,
but equality in our upper bound can only be attained for $2h=w$.

\subsection{Lower Bounds on the Geometric Dilation of Polygons}

In this subsection we apply the lower bound
(\ref{LowerDilationBoundDw-equ})
of Ebbers-Baumann et al.~\cite{egk-gdcpc-04ii}
to deduce lower bounds on the dilation of polygons with $n$ sides
(in special cases we proceed directly). We start with the case of a
triangle (and skip the easy proof):
\begin{lem} \label{LowerDilationBoundTriangle-lem}
For any triangle $C$, $\delta(C) \geq 2$. This bound is tight.
\end{lem}
\comment{
\begin{pf}
Denote by $a,b,c$ the sides of the triangle, and assume that
$a \leq b \leq c$.
Consider a halving chord of length $x$ which is
parallel to side $a$, that forms a triangle similar to $C$,
whose sides are $z$ (on $c$) and $y$ (on $b$).
We have
\[
  \frac{x}{a}=\frac{y}{b}=\frac{z}{c},
  {\rm \ \ and \ \ } a+c-z+b-y=z+y.
\]
Thus
\[
  x=\frac{a(a+b+c)}{2(b+c)}
  {\rm \ \ and \ \ } \delta(C) \geq \frac{a+b+c}{2x}=\frac{b+c}{a}
\geq 2 .
\]
It is easy to check that the bound holds with equality for an
equilateral triangle.
\hfill$\qed$
\end{pf}
}%
Equality is attained by equilateral triangles.
Note that plugging the inequality $D/w \geq 2/\sqrt{3}$
into (\ref{LowerDilationBoundDw-equ})
would only give
$\delta(C) \geq \pi/3+1/\sqrt{3} \approx 1.624$.
We continue with the case of centrally symmetric
convex polygons, for which we obtain a tight bound.
\begin{thm} \label{LowerDilationBoundSymmetricPolygon-thm}
If $C$ is a centrally symmetric convex $n$-gon (thus $n$ is even),
then
\[
 \delta(C) \geq \frac{n}{2} \tan{\frac{\pi}{n}} .
\]
This bound is tight.
\end{thm}
\begin{pf}
We adapt the proof of Theorem~22 in~\cite{egk-gdcpc-04ii}, which proves
inequality (\ref{LowerDilationBoundDw-equ}) for closed curves.
Since $C$ is centrally symmetric, it must contain a circle of radius~$r=h/2$.
It can easily be shown (using the convexity of the tangent function)
that the shortest $n$-gon containing such a circle is a regular $n$-gon.
Its length equals $2 r n \tan \pi/n$
which further implies that
\begin{equation}
  \label{LowerDilationBoundSymmetricPolygon-equ}
  \delta(C)
  \stackrel{\mbox{\scriptsize Lemma~\ref{DilationConvex-lem}}}{=}
  \frac{|C|}{2h}
  \geq
  \frac{hn\tan{\frac{\pi}{n}}}{2h}
  =
  \frac{n}{2} \tan{\frac{\pi}{n}} .
\end{equation}
The bound is tight for a
regular $n$-gon.
\hfill$\qed$
\end{pf}
In the last part of this section we address the case of arbitrary
(not necessarily convex) polygons.
Let $C$ be a polygon with
$n$ vertices, and let $C'=conv(C)$.
Clearly $C'$ has at most $n$ vertices.
By Lemma 9 in~\cite{egk-gdcpc-04ii},
$ \delta(C) \geq \delta(C')$.
Further on, consider
\[
  C''=\frac{C'+(-C')}{2}
  ,
\]
the convex curve obtained by \emph{central symmetrization} from $C'$
(see~\cite{yb-cf-61,egk-gdcpc-04ii}).
It is easy to check
that $C''$ is a convex polygon, whose number of vertices~$n''$ is even and at most
twice that of $C'$, therefore at most $2n$.
Because $\delta(C'') \leq \delta(C')$ by Lemma~16 in~\cite{egk-gdcpc-04ii},
we get
\[
  \delta(C) \geq \delta(C') \geq \delta(C'')
  \stackrel{\textrm{\scriptsize Theorem~\ref{LowerDilationBoundSymmetricPolygon-thm}},\;n''\leq 2n}{\geq}
  n \tan \frac{\pi}{2n}
\]
and obtain a lower bound
on the geometric dilation of any polygon with $n$ sides.
\begin{cor} \label{LowerDilationBoundPolygon-cor}
The geometric dilation of any polygon $C$ with $n$ sides satisfies
\[ \delta(C) \geq n \tan{\frac{\pi}{2n}} . \]
\end{cor}
This inequality does not seem to be tight, even for odd~$n$.
The dilation of a regular polygon~$C$ with an odd number~$n$ of vertices
can be calculated by using the fact that the curve~$C^*$ obtained by
the  halving-pair transformation is a regular $2n$-gon
whose dilation equals $\delta(C^*)=n \tan(\pi/2n)$.
Because the derivative (unit-)vectors $\dot{c}(t)$ and $\dot{c}(t+|C|/2)$
of an arc-length parameterization of~$C$ always enclose an angle $(n-1)\pi/n$,
we get
\[
  |C^*|
  =
  \int_0^{|C|} \frac{1}{2} \left(\dot{c}(t)-\dot{c}\left(t+\halfC\right)\right) \intd t
  =
  |C| \sin\left(\frac{n-1}{n} \frac{\pi}{2}\right)
  =
  |C| \cos \frac{\pi}{2n}
  .
\]
Because of $h(C^*)=h(C)$ and Lemma~\ref{DilationConvex-lem} we have $\delta(C)/\delta(C^*) = |C|/|C^*|$,
which results in
\[
  \delta(C)
  =
  \frac{1}{\cos \frac{\pi}{2n}} \delta(C^*)
  =
  \frac{1}{\cos\frac{\pi}{2n}} n \tan \frac{\pi}{2n}.
\]
This exceeds the lower bound of Corollary~\ref{LowerDilationBoundPolygon-cor}
by a factor of $1/\cos(\pi/2n) \approx 1+ \pi^2/(8n^2)$.

\section{Conclusion and Open Questions}
Our main result Theorem~\ref{ImprovedLowerBound-thm}
looks like a very minor improvement over the easier bound
$\Delta\ge \pi/2$, but it settles the question whether $\Delta > \pi/2$
and has required the introduction of new techniques.
Our approximations are not very far from optimal, and
we believe that new ideas are required to improve the lower bound to, say,
$\pi/2+0.01$.
An improvement of the constant $\Lambda=1.00001$ in the disk packing result
of~\cite{kkmv-atpe-99} (Theorem~\ref{DiskPacking-thm}) would of course immediately imply a
better bound for the dilation.
It should be emphasized here that Theorem~\ref{DiskPacking-thm}
holds for both finite and infinite packings.
As we use it only for finite packings,
it would be interesting to know if a substantially better result
could be obtained for this presumably easier case.

\comment{The authors of~\cite{kkmv-atpe-99} did not attempt
to optimize the parameters of their proof in order to get the strongest
possible bound, but rather tried to choose parameters which would permit
realistic pictures of the situations in the proof.
The choice of parameters in their proof is very delicate,
and it is not straightforward to improve it.
}

We do not know whether the link between disk packing
and dilation that we have established works in the
opposite direction as well: Can one construct a graph of small
dilation from a ``good'' circle packing (whose enlargement by a
``small'' factor covers a large area)?  If this were true
(in some meaningful sense which would have to be made precise)
it would mean
that a substantial improvement of the lower bound on dilation cannot
be obtained without proving, at the same time, a strengthening of
Theorem~\ref{DiskPacking-thm} with a larger constant than 1.00001.
\comment{
On the other hand, we know that also the upper bound
on $\Delta$ coming from the grid of Figure~\ref{DilationExample-fig}
can be improved, but we have not found any important improvement, yet.
}%
Overall, the gap between the lower bound $(1+10^{-11})\pi/2\approx 1.571$
and the upper bound $1.678$ remains a challenging problem.

As mentioned in Section~\ref{hAndA-sec},
we conjecture that the Rounded Triangle~$C_\triangle$ of
Figure~\ref{RoundedTriangleAndFlower-fig}a
is the convex curve of constant halving distance minimizing the area.


Finally, it would be nice to find a tight lower
bound on the geometric dilation of arbitrary (not only centrally-symmetric)
convex polygons.

{\em Acknowledgement.} We would like to thank
John Sullivan and Salvador Segura-Gomis for helpful discussions
and the anonymous referees for their valuable comments.

\addcontentsline{toc}{section}{References}

\appendix

\section{Proof of the Precise Bound in Lemma \ref{Ring-lem}}
\label{PreciseBound-sec}%
The main assumption of the lemma is
$\delta(C) \leq (1 + \eps)\pi/2$ for $\eps \leq 0.0001$.
Assume $H/h=(1+\beta)$.
Lemma~\ref{LowerDilationBoundC-lem}
yields the lower bound:
\[
  \delta(C)
  \geq
  \arcsin \frac{1}{1+\beta} + \sqrt{(1+\beta)^2-1}
  =
  \arcsin \frac{1}{ 1+\beta} + \sqrt{2\beta+\beta^2}
\]
We have $\beta \leq 0.01$, otherwise
this implies $\delta(C)> 1.0001\,\pi/2$, which contradicts the
assumption of the lemma.

It is well known that for $x \in [0,\pi/2]$,
\[
  \cos{x}
  \leq
  1-\frac{x^2}{2}+ \frac{x^4}{24}
  .
\]
By setting $x=\sqrt{2\beta}$,
we obtain the following inequality,
for the given $\beta$-range:
\[
  \sin{\left(\frac{\pi}{2}-\sqrt{2\beta}\right)}
  =
  \cos{\sqrt{2\beta}}
  \leq
  1-\beta+ \frac{\beta^2}{6}
  \stackrel{\scriptstyle \beta \leq 0.01}{\leq}
  1-\frac{\beta}{\beta+1}
  =
  \frac{1}{\beta+1}
  .
\]
Thus
\[
  \arcsin \frac{1}{1+\beta} \geq \frac{\pi}{2}-\sqrt{2\beta},
\]
and therefore
\begin{eqnarray*}
  \delta(C)
  & \geq &
  \frac{\pi}{2}-\sqrt{2\beta}+\sqrt{2\beta+\beta^2}
  =
  \frac{\pi}{2}+\frac{\beta^2}{\sqrt{2\beta}+\sqrt{2\beta+\beta^2}}
  \\ & \stackrel{\scriptstyle \beta \leq 0.01}{\geq} &
  \frac{\pi}{2}+\frac{\beta^2}{(\sqrt{2}+\sqrt{2.01})\sqrt{\beta}}
  \geq \frac{\pi}{2}+\frac{\beta^{3/2}}{3}.
\end{eqnarray*}
As a parenthesis, in our earlier estimate, equation (\ref{AsymptoticBound-equ}),
we have used only an asymptotic expansion for $\arcsin(1/(1+\beta))+\sqrt{(1+\beta)^2-1}$
without a precise bound on the error term.
Using the expansion in equation~(\ref{AsymptoticBound-equ}) one would probably
get a slightly better lower bound on $\Delta$,
but the improvement over $\pi/2$ would still be of the same order $10^{-11}$.

With our initial assumption, we get
\[
  \frac{\pi}{2}+\frac{\beta^{3/2}}{3}
  \leq
  \delta(C)
  \leq
  \frac{\pi}{2}(1+\eps),
\]
which yields
\begin{equation}
  \label{Beta-equ}
    \beta
    \leq
    \left(
      \frac{3\pi}{2}
    \right)^{2/3}
    \eps^{2/3}
    \leq
    2.9 \eps^{2/3}
    \stackrel{\scriptstyle \eps \leq 10^{-4}}{\leq}
    0.7 \eps^{1/2}.
\end{equation}
Lemma~\ref{BoundLengthM-lem} gives
\begin{equation}
  \label{BoundLengthM-equ}
  |M|
  \stackrel{\textrm{\scriptsize Lemma \ref{BoundLengthM-lem}}}{\leq}
   \frac{\pi h}{2} \sqrt{2\eps+\eps^2}
   \stackrel{\scriptstyle \eps \leq 10^{-4}}{\leq}
   2.24\,h\sqrt{\eps}
   .
\end{equation}
We have to bound the ratio~$R/r$ between the two concentric circles
containing~$C$.
\begin{eqnarray*}
  \frac{R}{r}
  & \;=\; &
  \frac{H/2 +|M|/4}{h/2 -|M|/4}
  \;=\;
  \frac{h(1+\beta) + |M|/2}{h - |M|/2}
  \\ & \;\stackrel{\textrm{\scriptsize \thetag{\ref{BoundLengthM-equ}}}}{\leq}\; &
  \frac{1+\beta + 1.12\sqrt{\eps}}{1 - 1.12\sqrt{\eps}}
  \;\stackrel{\mbox{\scriptsize \thetag{\ref{Beta-equ}}}}{\leq}\;
  \frac{1+ 1.82 \sqrt\eps }{1- 1.12 \sqrt\eps }
  \;\stackrel{\scriptstyle \eps \leq 10^{-4}}{\leq}\;
  1+3 \sqrt\eps
\end{eqnarray*}
This completes the proof of the precise bound in Lemma~\ref{Ring-lem}.
\hfill$\qed$

\section{Detailed Analysis of the Tightness Example for the Stability Result}
\label{DetailedTightnessExample-sec}%
In the end of Section~\ref{RingResult-sec} we defined
a cycle $C$ which shows that the coefficient~$3$
in Lemma~\ref{Ring-lem} cannot be smaller
than~$3/2$. Here, we discuss this in detail.

The norm of the derivative of $c$ in the given parameterization can be
calculated exactly:
\[
  |\dot{c}(\varphi)| =  \sqrt{1 + (64/9) s^2(1 - \cos^2(3 \varphi))}
= 1 + O(s^2)
\]
This means that the length of the curve piece $C[\varphi_1,\varphi_2]$
between two parameter values
$\varphi_1< \varphi_2$ is closely approximated by the difference of
parameter values.
\begin{equation}
  \label{LengthOK-equ}
  |C[\varphi_1,\varphi_2]|
 = \int_{\varphi_1}^{\varphi_2} |\dot{c}(\varphi)|\,\intd \varphi
= (\varphi_2-\varphi_1) (1 + O(s^2))
\end{equation}
In particular, the total length is
$2\pi + O(s^2)$.
It follows from \thetag{\ref{LengthOK-equ}}
that halving pairs are defined by parameter values
$\varphi$ and $\hat\varphi = \varphi \pm \pi \pm O(s^2)$.
The motion of $C$ can be decomposed into a circular orbit of the earth
and a local elliptic orbit of the moon:
\begin{eqnarray*}
  &&
  c(\varphi)
  =
  e(\varphi)+m(\varphi)
  ,\;\;
  e(\varphi)
  :=
  \binom{\cos \varphi}{\sin \varphi}
  ,
  \\ &&
  m(\varphi)
  :=
  s\cdot\left(\binom{\cos \varphi}{\sin \varphi}  \cos{3 \varphi}
  +
  \binom{-\sin \varphi}{\cos \varphi} (-\frac{1}{3} \sin {3 \varphi})\right)
  .
\end{eqnarray*}
Note that $m(\varphi)$ does not denote the midpoint curve but the moon's
curve with respect to the earth.
Points of ``opposite'' parameter values $\varphi$ and $\bar\varphi := \varphi + \pi$ have
exactly distance~$2$, since the terms in $m$ cancel:
$m(\varphi+\pi)=m(\varphi)$, and hence
\[
  c(\varphi)-c(\bar\varphi)
  =
  2 \binom{\cos \varphi}{\sin \varphi}.
\]
Halving distances can be estimated as follows:
\begin{eqnarray*}
  \lefteqn{|c(\varphi)-c(\hat\varphi)|^2
  =
  \left|
    \bigl(c(\varphi)-c(\bar\varphi)\bigr)
    +
    \bigl(c(\bar\varphi)-c(\hat\varphi)\bigr)
  \right|^2}
  \\ & = &
  |c(\varphi)-c(\bar\varphi)|^2
  +
  2
  \left\langle
    c(\varphi)-c(\bar\varphi),c(\bar\varphi)-c(\hat\varphi)
  \right\rangle
  +
  |c(\bar\varphi)-c(\hat\varphi)|^2
  \\ & = &
  4
  +
  2
  \left\langle
    \binom{\cos \varphi}{\sin \varphi},
    e(\bar\varphi)-e(\hat\varphi)
    +
    m(\bar\varphi)-m(\hat\varphi)
  \right\rangle
  +
  [O(s^2)(1+O(s^2))]^2
\end{eqnarray*}
The estimate for the last expression follows from
$|\bar\varphi-\hat\varphi|=O(s^2)$ and \thetag{\ref{LengthOK-equ}}.
The scalar product can be decomposed into two terms.
The first term can be evaluated directly:
\[
  \left\langle
    \binom{\cos \varphi}{\sin \varphi},
    e(\bar\varphi)-e(\hat\varphi)
  \right\rangle
  =
  -(1-\cos(\bar\varphi-\hat\varphi))
  =
  O(\bar\varphi-\hat\varphi)^2
  =
  O(s^4)
\]
The second term can be bounded by noting that the moon's speed is bounded:
$|\dot{m}|=O(s)$.
\[
  \left\langle
    \binom{\cos \varphi}{\sin \varphi},
    m(\bar\varphi)-m(\hat\varphi)
  \right\rangle
  \le
  1 \cdot |m(\bar\varphi)-m(\hat\varphi)|
  =
  O(s) \cdot |\bar\varphi-\hat\varphi|
  =
  O(s^3)
\]
Putting everything together, every halving distance is bounded as follows:
\[
  |c(\varphi)-c(\hat\varphi)|
  =
  \sqrt{4- O(s^4) \pm O(s^3)+ O(s^4)}
  =
  2 \pm O(s^3)
  .
\]
$H$ and~$h$ are bounded by the same estimate.

A more precise estimate for the length is
$|C|=2\pi(1+ \frac{16}9s^2 + O(s^4))$.
Substituting this into the derivation
at the end of Section~\ref{RingResult-sec} gives a dilation of
$\delta(C) = (1+\eps)\pi/2$ with $\eps=\frac{16}9s^2 + O(s^3)$.
The ratio of the radii of the enclosing ring is
$(1+s)/(1-s)=1+2s+O(s^2)=1+\frac{3}{2}\sqrt\eps+O(\eps)$.
This means that the coefficient $3$ of $\sqrt\eps$
in Lemma~\ref{Ring-lem} cannot be reduced below $3/2$.
\hfill$\qed$

\section{Parameterization of the Rounded Triangle}%
\label{ParameterizationRoundedTriangle-sec}%
In this section we solve a differential equation to give
a parameterization $c(t)=(x(t),y(t))$ of a part of~$C_\triangle$.
It is the first half rounded piece.
As the curve $C_\triangle$ contains three curved pieces, it consists of
six halves like the one described in the following.
Together with the straight line segments of the same length 
they build the whole Rounded Triangle.
The piece examined here starts at $c(0):=(0,-0.5)$ and it is determined by the
two conditions
\begin{eqnarray*}
  \sqrt{(x(t)-t)^2+(y(t)-0.5)^2}
  \equiv
  1
  & \;\;\mbox{and}\;\; &
  \sqrt
  {
    x'(t)^2
    +
    y'(t)^2
  }
  \equiv
  1
  \;.
\end{eqnarray*}
Solving the first one for $y$ and taking the derivative
with respect to $t$, and solving the second
one for $\pseudofrac{\intd y}{\intd t}$ yields
\[
  -\frac
  {(x(t)-t)\left(x'(t) - 1\right)}
  {\sqrt{1-(x(t)-t)^2}}
  =
  y'(t)
  =
  \sqrt{1-x'(t)^2}
  \;.
\]
By squaring we get the quadratic equation
\[
  x'(t)^2
  -
  2(x(t)-t)^2x'(t)
  +
  2(x(t)-t)^2-1
  =
  0
  \;.
\]
As the second possible solution $x'(t) \equiv 1$
does not make any sense in this context, we get
\[
  x'(t)
  =
  2(x(t)-t)^2-1
  \;.
\]
This differential equation
with the constraint $x(0)=0$ yields
\begin{equation}
  \label{eq:parame}
  x(t)
  =
  t
  -
  \frac
  {e^{4t}-1}
  {e^{4t}+1}
  \;\;\mbox{and}\;\;
  y(t)
  =
  -2\frac{e^{2t}}{e^{4t}+1}+0.5\;.
  \end{equation}
The first of the twelve pieces ends when
the tangent has reached an angle of $30^\circ$
with the $y$-axis, i.e. $\sin(x'(t_1))=\pi/6$,
$x'(t_1)=-1/2$.
Using the formula for $x'(t_1)$ and
substituting $z := e^{4t_1}$, we get
\[
  \frac{z^2-6z+1}{(z+1)^2} = -\frac{1}{2}
\]
which has the solution $z=(5/3)\pm\sqrt{(5/3)^2-1}=(5/3)\pm(4/3)$.
As we are looking for a positive solution, we get
$t_1 = \ln3/4$.
The whole closed curve consists of twelve parts of this length.
Hence, its perimeter and dilation are given by
\[
  \abs{C_\triangle} = 12 \frac{\ln 3}{4} = 3\ln 3 \approx 3.2958
  \;\;,\;\;\;\;\;\;
  \delta\left(C_\triangle\right) = \frac{3}{2}\ln 3 \approx 1.6479
  \;.
\]
\end{document}